**Solutions to Two Open Problems in Geometric Group Theory**

**by**

**Jordan A. Sahattchieve**

**A dissertation submitted in partial fulfillment
of the requirements for the degree of
Doctor of Philosophy
(Mathematics)
in The University of Michigan
2012**

**Doctoral Committee:**

       **Professor G. Peter Scott, Chair
       Professor Richard D. Canary
       Professor Stephen M. DeBacker
       Professor Karen E. Smith
       Associate Professor James P. Tappenden**

To my parents

THE CONTENT OF THIS PAGE HAS BEEN REMOVED

# Prologue

I learned about bounded packing in polycyclic groups from a talk Chris Hruska gave at a special session of the American Mathematical Society at the University of Illinois in Urbana-Champaign. At the time, I had gotten stuck on a rather hard problem which involved understanding the convex hull in arbitrary CAT(0) spaces. I had a simple idea of constructing the convex hull of a set in steps by adjoining geodesic segments until a convex set is reached. I realized that if one is able to show that the number of steps needed to produce the convex hull of a quasiconvex subset is finite, one obtains a bound on the maximal distance between a point in the convex hull of the set and the nearest point in the set itself, which implies cocompactness of certain group actions thus settling a special case of a longstanding open problem in the theory of CAT(0) groups. Later, I realized that Hermann Brünn had already investigated this process and had given upper bounds on the number of steps required for the process to terminate in certain Banach spaces. I was unable to procure his original paper in German, which is why I proved a result similar to, but possibly weaker than his own for Euclidean spaces.

At the time I was ready to announce my solution to the bounded packing problem in polycyclic groups, I received communication that the problem had been solved.

# Table of Contents

# ABSTRACT

## Solutions to Two Open Problems in Geometric Group Theory

### by

### Jordan A. Sahattchieve

## Chair: G. Peter Scott

We introduce a method for analyzing the convex hull of a set in non-positively curved piecewise Euclidean polygonal complexes and we apply this method to prove that with the usual action of $F_m \times \mathbb{Z}^n$ on the metric product of a tree with $\mathbb{R}^n$, every quasiconvex subgroup of $F_m \times \mathbb{Z}^n$ is convex. This answers the question whether a quasiconvex subgroup of a CAT(0) group is a CAT(0) group in the affirmative for the groups $F_m \times \mathbb{Z}^n$. We also prove bounded packing in a special class of polycyclic groups, and we introduce the notion of coset growth and provide a bound for the coset growth of uniform lattices in *Sol*.

# Chapter 1

# Introduction and Synopsis

**Original Results**

This dissertation contains the solutions to two unrelated open problems in geometric group theory.

The first part of this dissertation contains a solution to a special case of the bounded packing problem for polycyclic groups. Without going into too much detail, having the property of bounded packing for a group means that the cosets of any of its subgroups are spread out, meaning that large enough collections of distinct cosets are guaranteed to contain a pair of cosets which are "uniformly" far from each other. We use the technique introduced in Section 3.5 for the purpose of showing bounded packing in certain extensions of $\mathbb{Z}^n$ to find a bound for the coset growth in the non-nilpotent groups $\mathbb{Z}^2 \rtimes \mathbb{Z}$ in Section 3.7.2.

The second part of this dissertation explores the notion of quasiconvexity in CAT(0) groups, particularly in the groups $F_m \times \mathbb{Z}^n$. This work was motivated by the question of whether subgroups of CAT(0) groups are themselves CAT(0) groups. While the answer to this question is an emphatic *no*, one can still hope that answer may become "yes" after attaching enough adjectives to the group and the subgroup in question. We show that in the case of the groups $F_m \times \mathbb{Z}^n$ acting in the usual way on Tree×$\mathbb{R}^n$, this is indeed the case. In fact, we show something stronger - we show that any *quasiconvex* subgroup of $F_m \times \mathbb{Z}^n$ acts on the convex hull of any of its orbits with a compact quotient. In more rigorous terms, our result says that for the groups $F_m \times \mathbb{Z}^n$, *quasiconvex* $\Rightarrow$ *convex*.

## Synopsis

The main body of this dissertation contains three chapters.

Chapter 2 is an introduction to the following notions in geometric group theory which we make use of in the course of proving our results: Section 2.1 provides a very quick and informal introduction to the world of metric geometry and CAT(0) groups. Section 2.2 goes further and provides a relatively self-contained introduction to CAT(0) metric spaces and their isometry groups, and $\delta$-hyperbolic metric spaces; much of the the material in this section can be found in [3]. Section 2.3 is an introduction to the notions of ends of groups and ends of group pairs, and also contains a summary of the basic properties of Sageev's well-known cubing construction.

Chapter 3 contains a solution to the bounded packing problem in certain extension of $\mathbb{Z}^n$ by $\mathbb{Z}$ and introduces the notion of growth of group pairs. Section 3.1 discusses the origins of the bounded packing property, while Sections 3.2 and 3.3 provide rudimentary background on bounded packing and polycyclic groups respectively. Following a discussion of the dynamics of linear maps on $\mathbb{R}^n$ in Section 3.4, the chapter proceeds to analyze the diagonal action of $\mathbb{Z}$ on $\mathbb{Z}^n$ in Section 3.5. The chapter then concludes with a bound on the growth of the group pair $(\mathbb{Z}^2 \rtimes \mathbb{Z}, \langle t \rangle)$ found in Section 3.7.

In Chapter 4, we provide an (affirmative) answer the question whether a quasi-convex subgroup of a CAT(0) group is convex for the groups $F_m \times \mathbb{Z}^n$ with respect to the usual action of these groups on Tree$\times \mathbb{R}^n$. Section 4.1 contains a discussion of the notion of quasiconvexity in CAT(0) spaces and CAT(0) groups while Section 4.2 outlines the iterative construction of the convex hull in CAT(0) spaces, which we use in the following Section 4.3 to prove convexity of all quasiconvex subgroups of $F_m \times \mathbb{Z}^n$. Section 4.3 also demonstrates that all quasiconvex subgroups of $F_m \times \mathbb{Z}^n$ are virtually of the form $A \times B$ for $A \leq F_m$ and $B \leq \mathbb{Z}^n$.

All citations made throughout this dissertation can be found in the bibliography section at the end of the document.

# Chapter 2

# Background

## 2.1 Non-positive Curvature in Metric Spaces

One would like to study the group of isometries of a metric space in a manner analogous to the study of discrete groups of isometries of Euclidean and hyperbolic spaces. While for arbitrary metric spaces this is perhaps a futile endeavor, the idea to study such isometry groups has proven to be a fertile field of study for geodesic metric spaces which are non-positively curved in a certain sense which we will define shortly.

The notion of non-positive curvature has its roots in the realm of differential geometry. While there is no ambiguity in what it means for a Riemannian manifold to be non-positively curved, it is not clear how this notion is to be translated to arbitrary geodesic metric spaces. In fact, there are two distinct, yet intimately related, notions of non-positive curvature.

One way to approach the issue is by analogy with Hadamard manifolds (manifolds of non-positive sectional curvature). It has long been known that geodesic triangles in a Hadamard manifold whose sectional curvature is bounded above by $\kappa \leq 0$ everywhere are thinner in comparison to triangles in the space form of constant curvature $\kappa$ (which we shall denote by $\mathcal{M}^\kappa$). Let $x, y$ and $z$ be points on our favorite Hadamard manifold whose curvature is bounded above by $\kappa$. We can then find three points $\overline{x}, \overline{y}, \overline{z} \in \mathcal{M}^\kappa$, such that their relative distances are the same as the relative distances between $x, y$ and $z$, more precisely, $d(x, y) = d(\overline{x}, \overline{y})$, $d(x, z) = d(\overline{x}, \overline{z})$, and $d(y, z) = d(\overline{y}, \overline{z})$. Then, we have the following classical comparison inequality: if $a$ and $b$ are any two points on the geodesic triangle $\Delta(x, y, z)$, and $\overline{a}, \overline{b}$ are the corresponding points on the $\mathcal{M}^\kappa$ comparison triangle $\Delta(\overline{x}, \overline{y}, \overline{z})$, we have $d(a, b) \leq d(\overline{a}, \overline{b})$. In the differential geometric setting, one proves this as a consequence of the non-positive sectional curvature condition. In an arbitrary geodesic metric space, where there is no notion of curvature, we can use this comparison inequality as a means of detecting the salient features of non-positive curvature. We define a *CAT(0) space* to be a geodesic metric space in which geodesic triangles are thinner than the corresponding triangles in $\mathcal{M}^0 = \mathbb{R}^2$. The letters in the acronym $CAT(0)$ stand for the last initials of the names of the mathematicians Cartan, Alexandrov, and Toponogov.

One can similarly define metric spaces of curvature bounded above by $\kappa$ for arbitrary $\kappa \in \mathbb{R}$. This definition is due to A. D. Alexandrov, who called them $\mathcal{R}_\kappa$ domains until Gromov renamed them CAT($\kappa$) spaces. We now define a *CAT(0) group $G$* to be a group that acts properly and cocompactly by isometries on a CAT(0) space.

Another way to define negative curvature in a metric space is due to Gromov via the $\delta$-slim triangles condition. A group $\Gamma$ is said to be Gromov hyperbolic if its Cayley graph is a $\delta$-hyperbolic metric space, for some $\delta$, meaning that given any geodesic triangle $\Delta(x, y, z)$, each side is contained in the $\delta$-neighborhood of the union of the other two sides.

Finally, there is a more combinatorial approach to geometric group theory described in the seminal paper by Scott and Wall [31] which is based on covering space theory and analyzes the structure of finitely generated groups by examining their actions on trees.

## 2.2 CAT(0) and $\delta$-Hyperbolic Groups

### 2.2.1 Basic Definitions

In this section we define the basic notion of a metric, metric space, and isometry, and present some fundamental motivating examples for the study of CAT(0) spaces.

Let $X$ be a set. A function $X \times X \to \mathbb{R}$ is called a *pseudometric* if for any $x, y, z \in X$:

1. $d(x, y) \geq 0$

2. $d(x, y) = d(y, x)$

3. $d(x, z) \leq d(x, y) + d(y, z)$

If in addition to the above, $d$ also satisfies $d(x, y) \neq 0$ whenever $x \neq y$, we call $d$ a *metric*, and the pair $(X, d)$ a *metric space*. Recall from basic topology that the metric $d$ induces a natural topology on $X$, which has as basis the open balls $B(x, r) = \{y \in X : d(x, y) < r\}$, see for example [20]. This topology makes the metric space $X$ into a *perfectly normal space*, in other words, a topological space in which any two disjoint non-empty closed sets can be separated by a continuous function $f : X \to \mathbb{R}$.

We call the set $\overline{B}(x, r) = \{y \in X : d(x, y) \leq r\}$ the *closed ball* of radius $r$ centered at $x$ and note that $\overline{B}(x, r)$ may be strictly larger than the closure of the open ball $B(x, r)$. A metric space which has the property that $\overline{B}(x, r)$ is compact for every $x \in X$ and every $r > 0$ is called a *proper* metric space. The first examples of metric spaces one encounters in introductory courses on analysis, namely the finite dimensional Euclidean spaces $\mathbb{R}^n$, are all proper. So are also all Riemannian manifolds. For examples of spaces which are not proper, one usually looks at infinite dimensional Banach spaces. Since Banach spaces are an interesting class of spaces on which we can do geometry, let us recall their definition here.

**Definition 2.2.1.** Let $V$ be a vector space over $\mathbb{R}$. A *norm* on $V$ is a function $\|\cdot\| : V \to \mathbb{R}_{\geq 0}$, satisfying:

1. $\|a\mathbf{v}\| = |a| \, \|\mathbf{v}\|$, for all $a \in \mathbb{R}$, $\mathbf{v} \in V$

2. $\|\mathbf{v} + \mathbf{w}\| \leq \|\mathbf{v}\| + \|\mathbf{w}\|$, for all $\mathbf{v}, \mathbf{w} \in V$

3. If $\|\mathbf{v}\| = 0$, then $\mathbf{v} = \mathbf{0}$

Given a normed linear space $V$, one obtains a metric on $V$ by setting $d(\mathbf{v}, \mathbf{w}) = \|\mathbf{v} - \mathbf{w}\|$. A normed linear space $V$ which is complete in the metric $d$ induced by its norm is called a *Banach space*. For basics on Banach spaces, see [18].

Recall that $l^p(\mathbb{R})$ is the Banach space of all functions $\mathbf{x} : \mathbb{N} \to \mathbb{R}$, such that $\|\mathbf{x}\|_p = (\Sigma_i |x_i|^p)^{1/p} < \infty$. The spaces $L^p(\mathbb{R}^n)$ are the normed linear spaces of Lebesgue integrable real-valued functions on $\mathbb{R}^n$ with the norm $\|f\|_p = \left(\int |f|^p \, d\mu\right)^{1/p} < \infty$. The integral is taken over all of $\mathbb{R}^n$ with respect to the Lebesgue measure $d\mu$. For basic facts on the $L^p$ and $l^p$ spaces such as completeness and duality properties, see the timeless classic by Walter Rudin [25]. By contrast with the finite dimensional vector spaces $\mathbb{R}^n$, none of the infinite dimensional Banach spaces $l^p(\mathbb{R})$, $L^p(\mathbb{R})$ are proper, see [18].

Having defined the notion of a metric, it is only natural to consider the maps between metric spaces which preserve distances.

**Definition 2.2.2.** Let $(X_1, d_1)$ and $(X_2, d_2)$ be two metric spaces. A map $f : X_1 \to X_2$ is called an *isometric embedding*, if $d_1(x, y) = d_2(f(x), f(y))$, for all $x, y \in X_1$.

Note that an isometric embedding is always a continuous injection. If in addition to being an isometric embedding, $f$ is also onto, we call $f$ an *isometry*. We note that the isometries of a metric space $X$ to itself form a group, called the isometry group of $X$, denoted by $Isom(X)$. It is this group that is of interest to us. As we mentioned in Section 2.1, the group $Isom(X)$ is hard to work with and can often turn out to be trivial, in other words, the space $X$ may have very few symmetries. Usually, one is interested in the isometry group of spaces which have more structure than just a metric structure. For example, a very active field of research is the study of *Kleinian groups*, or discrete groups of isometries of $\mathbb{H}^3$. In addition to having a smooth structure $\mathbb{H}^3$ is also simply connected and has constant curvature which equals $-1$. Kleinian groups are an important example of both CAT(0) and hyperbolic groups, which we define below, and provide a source of test questions for the theory of negatively curved groups. Let us work toward describing the most general sort of spaces on which a sufficiently rich geometric structure exists. Intuitively, in order to do geometry, one needs lines and angles. An isometric embedding $\lambda$ of an interval $I \subseteq \mathbb{R}$ into the metric space $X$ is called a *geodesic*. For us, geodesics will be the equivalent of lines in Euclidean geometry. Often, one blurs the distinction between the map $\lambda$ and its image. If $x, y \in X$, one calls the image of a geodesic $\lambda : [a, b] \to X$ such that $\lambda(a) = x$ and $\lambda(b) = y$ a *geodesic segment* between $x$ and $y$ and denotes it by $[a, b]$ whenever there is no ambiguity in doing so. Given a geodesic $\lambda : I \to \mathbb{R}$, one calls $\lambda$ a *linearly parametrized geodesic* or geodesic *parametrized proportional to arc length* if there exists a constant $c$ such that $d(\lambda(t_1), \lambda(t_2)) = c|t_1 - t_2|$, for all $t_1, t_2 \in I$. A *geodesic ray* is a map $\lambda : [0, \infty) \to X$, such that $d(\lambda(t_1), \lambda(t_2)) = |t_1 - t_2|$, for all $t_1, t_2 \geq 0$. In other words, a geodesic ray is an isometric embedding of the positive real line into $X$ parametrized with unit speed. A *local geodesic* in $X$ is a map $\lambda : I \to X$ having the property that for every $t \in I$ there exists $\epsilon > 0$ such that $d(\lambda(t_1), \lambda(t_2)) = |t_1 - t_2|$, for all $t_1, t_2 \in I$ with $|t - t_1| + |t - t_2| \leq \epsilon$.

**Definition 2.2.3.** A metric space $(X, d)$ is called a *geodesic metric space*, or simply a *geodesic space*, if every two points in $X$ can be joined by a geodesic. The geodesic metric space $(X, d)$ is called *uniquely geodesic* if given $x, y \in X$ there is only one geodesic joining $x$ to $y$.

**Example 2.2.1.** The space $\mathbb{R}^n$ with the Euclidean or $l^2$ metric is a uniquely geodesic metric space. The map $[0,1] \to \mathbb{R}^n$ defined by $\lambda(t) = t\mathbf{x} + (1-t)\mathbf{y}$ parametrizes the unique geodesic segment joining the points $\mathbf{x}, \mathbf{y} \in \mathbb{R}^n$. Henceforth, we shall denote $\mathbb{R}^n$ with the metric coming from the $l^2$ norm by $\mathbb{E}^n$. The hyperbolic spaces $\mathbb{H}^n$ are also uniquely geodesic.

Perhaps the most familiar examples of geodesic metric spaces which are not uniquely geodesic are the spheres $(\mathbb{S}^n, d)$. In order to avoid a discussion of Riemannian structures here, we shall simply define the space $\mathbb{S}^n$ to be the unit sphere in $\mathbb{R}^n$ with the angular metric: if we denote the origin by $O$, given any $A, B \in \mathbb{S}^n$, we set $d(A, B)$ to be the acute Euclidean angle between the line segments $OA$ and $OB$. If $OA \perp OB$, we define $d(A, B) = \frac{\pi}{2}$. It is an easy exercise in spherical geometry to show that given any two points $A, B \in \mathbb{S}^n$ with $d(A, B) < \frac{\pi}{2}$, there exists a unique geodesic segment joining $A$ to $B$, namely, the shorter part of the great circle passing through these two points. On the other hand, if $d(A, B) = \frac{\pi}{2}$, there are infinitely many geodesic segments joining $A$ to $B$.

Recall that the space $l_n^\infty(\mathbb{R})$ is defined to be $\mathbb{R}^n$ endowed with the $l^\infty$ norm: $\|\mathbf{x}\|_\infty = \max\{x_1, ..., x_n\}$. It is an instructive exercise to show that $\|\mathbf{x}\|_p \to \|\mathbf{x}\|_\infty$, for any $\mathbf{x} \in \mathbb{R}^n$, as $p \to \infty$. Here is a more interesting example of a geodesic metric space which is not uniquely geodesic:

**Example 2.2.2.** The graph of any monotone function $I \to \mathbb{R}$ is a geodesic segment in $\mathbb{R}^2$ for the $l^\infty$ metric. In any Banach space $V$, the affine linear segments $t \to t\mathbf{x} + (1-t)\mathbf{y}$, for $\mathbf{x}, \mathbf{y} \in V$, $0 \le t \le 1$ are always geodesics (see [21], Proposition 5.3.7) and are called *affine geodesics*.

Whenever $(X, d)$ is a uniquely geodesic space, we shall say that a subset $C \subseteq X$ is *geodesically convex*, or simply *convex*, if for every $x, y \in C$, the geodesic segment joining $x$ to $y$ is completely contained in $C$. If $(X, d)$ is a geodesic space which is not uniquely geodesic, then the definition of convexity becomes slightly more complicated.

A subset $C \subseteq X$ is called *completely geodesically convex* if given $x, y \in C$, **every** geodesic segment joining $x$ to $y$ is contained in $C$. In Example 2.2.2, the affine geodesic segment $[x, y]$ joining any two points in $l_2^\infty(\mathbb{R})$ is an example of a convex set which is not completely geodesically convex.

Given a subset $S \subseteq X$, the smallest convex subset of $X$ which contains $S$ is called the *convex hull* of $S$. We shall denote this convex set by $conv(S)$. It is clear from the definition that $conv(S)$ is the intersection of all convex subsets of $X$ which contain $S$. In what follows, it will be necessary to consider subsets which although not convex themselves are not far from being convex. We call a subset $S \subseteq X$ *quasiconvex* if the convex hull of $S$ is contained in a bounded neighborhood of $S$, symbolically $conv(S) \subseteq \mathcal{N}_\nu(S)$, for some $\nu > 0$.

**Example 2.2.3.** (The Hilbert geometry) In this example, we recall an important construction due in part to Hilbert (1895) as a generalization of the Klein model of hyperbolic space. Most of the material here can be found in [21].

Recall that given four ordered points $a_1, a_2, b_1, b_2$ in $\mathbb{R}^n$ satisfying $a_1 \neq b_1$ and $a_2 \neq b_2$, one defines their *cross ratio* $[a_1, a_2, b_1, b_2]$ by the formula

$$[a_1, a_2, b_1, b_2] = \frac{d(a_2, b_1)d(a_1, b_2)}{d(a_1, b_1)d(a_2, b_2)}$$

where $d$ is the standard Euclidean metric on $\mathbb{R}^n$. Now, let $A$ be a nonempty bounded open convex set in $\mathbb{R}^n$. Here, "convex" means convex relative to the Euclidean metric. Given $a_1, a_2 \in A$, we consider the Euclidean line $L$ containing them. This line intersects the topological boundary $\partial A$ of $A$ in exactly two points. Let us denote these points of intersection by $b_1$ and $b_2$, with the labeling chosen in a way so that $b_1, a_1, a_2, b_2$ are aligned in that order on $L$. We define the *Hilbert metric* $h_A$ on $A$ by

$$h_A(a_1, a_2) = \begin{cases} \ln[a_1, a_2, b_1, b_2] & \text{if } a_1 \neq a_2, \\ 0 & \text{if } a_1 = a_2. \end{cases}$$

One then proves that:

1. The map $h_A : A \times A \to \mathbb{R}_{\geq 0}$ is a metric.

2. The metric space $(A, h_A)$ is a proper geodesic metric space.

3. Each affine segment in $A$ is a geodesic segment for $h_A$.

4. The metric space $(A, h_A)$ is uniquely geodesic if and only if $\partial A$ does not contain a pair of affine segments that span a two-dimensional affine plane.

The Hilbert geometries provide us with a family of examples of proper geodesic metric spaces which are not uniquely geodesic and in which each affine segment is geodesic. This family of examples is largely distinct from the family of linear spaces possessing the same properties to which Example 2.2.2 belongs.

### 2.2.2 Comparison Triangles, Angles, and Rectifiable Curves

The material in this section follows closely the account in [3].

Let $X$ be a metric space, and let $a, b, c \in X$ be three points. Because the distances $d(a, b)$, $d(b, c)$, and $d(a, c)$ satisfy the triangle inequality in $X$, we can find three points $\overline{a}, \overline{b}, \overline{c} \in \mathbb{E}^2$ having the property that $d(a, b) = d(\overline{a}, \overline{b})$, $d(b, c) = d(\overline{b}, \overline{c})$, and $d(a, c) = d(\overline{a}, \overline{c})$. It is an easy exercise in Euclidean geometry to show that the triangle $\overline{\Delta}(\overline{a}, \overline{b}, \overline{c})$ is unique up to an isometry. If additionally $X$ is a uniquely geodesic metric space, we can join $a$, $b$, and $c$ by geodesic segments to obtain the *geodesic triangle* $\Delta(a, b, c) = [a, b] \cup [b, c] \cup [a, c]$. In this case, we say that $\overline{\Delta}(\overline{a}, \overline{b}, \overline{c})$ is a *comparison triangle* for $\Delta(a, b, c)$. Regardless of whether $X$ is geodesic or not, we shall define the *comparison angle* between $b$ and $c$ at $a$, denoted by $\overline{\angle}_a(b, c)$, to be the interior angle of $\overline{\Delta}(\overline{a}, \overline{b}, \overline{c})$ at $\overline{a}$. Now, we are ready to state the definition of angle:

**Definition 2.2.4.** (Alexandrov angle) Let $X$ be a metric space (not necessarily geodesic), and let $\lambda_1, \lambda_2 : [0, L] \to X$ be two geodesic segments starting at a common point $a \in X$, i.e. $\lambda_1(0) = \lambda_2(0) = a$. Given $t_1, t_2 \in [0, L]$, we consider the comparison angle $\overline{\angle}_a(\lambda_1(t_1), \lambda_2(t_2))$. The *Alexandrov angle*, or simply the angle between the geodesic segments $\lambda_1$ and $\lambda_2$, denoted by $\angle(\lambda_1, \lambda_2)$ is defined by the formula: $\angle(\lambda_1, \lambda_2) = \limsup_{t_1, t_2 \to 0} \overline{\angle}_a(\lambda_1(t_1), \lambda_2(t_2))$.

Using the law of cosines on the comparison triangle $\overline{\Delta}(\overline{a}, \overline{\lambda_1}(t_1), \overline{\lambda_2}(t_2))$, we obtain

$$\angle(\lambda_1, \lambda_2) = \lim_{t_1, t_2 \to 0} \cos^{-1}\left(\frac{t_1^2 + t_2^2 - [d(\lambda_1(t_1), \lambda_2(t_2))]^2}{2t_1 t_2}\right) \tag{2.1}$$

We should issue a word of caution regarding the Alexandrov angle. While it is true that the Alexandrov angle between two geodesic segments which are in the same germ will equal 0, the converse is not true, see [3]. Recall that two geodesic segments $\lambda_1, \lambda_1 : [0, L] \to X$ are said to be in the same germ if $\lambda_1(t) = \lambda_2(t)$ for $t < \epsilon$, for some $\epsilon > 0$. The relation of being in the same germ is an equivalence relation on the set of geodesic segments issuing at a point.

Another central notion in metric geometry is the notion of the length of a path:

**Definition 2.2.5.** (Length) Let $X$ be a metric space. A continuous map from a compact interval $I \subset \mathbb{R}$ to $X$ is called a *path* in $X$. We define the *length* of $\gamma$ to be the positive number $l(\gamma) = \sup \sum_{i=0}^{n-1} d(\gamma(t_i), \gamma(t_{i+1}))$, where the supremum is taken over all possible partitions $t_0 \leq t_1 \leq \cdots \leq t_n$ of $I$ (with no bound on $n$). Here, $t_0$ and $t_n$ are the left and the right endpoints of $I$.

Even though $\gamma$ is a continuous map from a compact interval, $\gamma$ may not have finite length. Whenever $l(\gamma) < \infty$, we shall say that $\gamma$ is *rectifiable*.

### 2.2.3   CAT (0) and $\delta$-Hyperbolic Spaces

We are now ready to define one of our central objects of study.

**Definition 2.2.6.** (CAT(0) Space) Let $(X, d)$ be a metric space. Let $\Delta(a, b, c)$ be a geodesic triangle in $X$, and let $\overline{\Delta}(\overline{a}, \overline{b}, \overline{c})$ be a comparison triangle. We shall say that $\Delta(a, b, c)$ satisfies the *CAT(0) comparison inequality* if given any $p, q \in \Delta(a, b, c)$ we have $d(p, q) \leq d(\overline{p}, \overline{q})$, where $\overline{p}, \overline{q} \in \overline{\Delta}(\overline{a}, \overline{b}, \overline{c})$ are the corresponding comparison points uniquely determined by the conditions $d(a, p) = d(\overline{a}, \overline{p})$, and $d(a, q) = d(\overline{a}, \overline{q})$. We say that the metric space $(X, d)$ is a *CAT(0) space* if every geodesic triangle in $X$ satisfies the CAT(0) comparison inequality.

It is immediate from the definitions that a geodesic CAT(0) space is uniquely geodesic. Let us state some basic properties of CAT(0) spaces.

Recall that a function $f : I \to \mathbb{R}$ defined on an interval $I \subseteq \mathbb{R}$ is called *convex* if the graph of $f$ lies below the line segment joining any two points on it. More precisely, $f((\tau-1)t_1 + \tau t_2) \leq (\tau-1)f(t_1) + \tau f(t_2)$, for any $t_1, t_2 \in I$, and $\tau \in [0, 1]$. Equivalently, the function $f$ is convex if its *epigraph* $Ep(f) = \{(x, y) \in I \times \mathbb{R} : f(x) \leq y\}$ is a convex subset of $\mathbb{E}^2$, see [21]. We shall say that a function $f : X \to \mathbb{R}$ defined on a geodesic metric space is *convex* if for any geodesic segment $\lambda : I \to X$ parametrized proportional to arc length, the function $f \circ \lambda : I \to \mathbb{R}$ is convex.

We are now prepared to state:

**Proposition 2.2.4.** *If $(X, d)$ is a CAT(0) metric space, then the metric $d : X \times X \to \mathbb{R}_{\geq 0}$ is convex.*

*Proof.* For a complete proof we refer the reader to Proposition 2.2 in [3]. However, this result is essential for the work in Chapter 4, so we give a proof of the following special case: Let $\lambda_1, \lambda_2 : [0, 1] \to X$ be two geodesics parametrized proportional to arc length such that $\lambda_1(0) = \lambda_2(0)$. Consider a comparison triangle $\overline{\Delta} \in \mathbb{E}^2$ for $\Delta(\lambda_1(0), \lambda_1(1), \lambda_2(1))$. For $0 \leq t \leq 1$, we have $d(\overline{\lambda_1(t)}, \overline{\lambda_2(t)}) = td(\overline{\lambda_1(1)}, \overline{\lambda_2(1)}) = td(\lambda_1(1), \lambda_2(1))$. The CAT(0) inequality gives $d(\lambda_1(t), \lambda_2(t)) \leq d(\overline{\lambda_1(t)}, \overline{\lambda_2(t)})$, hence $d(\lambda_1(t), \lambda_2(t)) \leq td(\lambda_1(1), \lambda_2(1))$. Since $d(\lambda_1(0), \lambda_2(0)) = 0$, this concludes the proof of the special case. $\qquad \square$

Another fundamental property of CAT(0) spaces is the existence of orthogonal projections:

**Proposition 2.2.5.** *Let $X$ be a CAT(0) space, and let $C \subseteq X$ be a convex subset, complete in the induced metric. Then, for every $x \in X$, there exists a unique point $pr_C(x) \in C$ such that $d(x, pr_C(x)) = \inf_{y \in C} d(x, y)$. Further, the map $pr_C : X \to C$ is a distance non-increasing retraction of $X$ onto $C$.*

*Proof.* See Proposition 2.4 in [3]. □

An important consequence of Proposition 2.2.5 is the following:

**Corollary 2.2.6.** *The map $h : X \times [0, 1] \to X$ which sends $(x, t)$ to the point at distance $td(x, pr_C(x))$ from $x$ along the geodesic segment $[x, pr_C(x)]$ is a continuous homotopy from $id_X : X \to X$ to $pr_C : X \to C$. In particular, every CAT(0) space is contractible.*

*Proof.* See Proposition 2.4 in [3]. □

Having defined what a CAT(0) space is, it is natural to look for examples of CAT(0) spaces different from the Euclidean spaces $\mathbb{E}^n$. Perhaps the most historically important such are furnished by Cartan's celebrated theorem which we briefly recall. While Cartan's result is more general than what we are able to state here, we make no use of it and will therefore refer the interested reader to Cartan's original paper [6] or the exposition in Bridson and Haefliger [3].

**Hadamard Manifolds**

Let $M$ be a smooth manifold together with a correspondence $p \to \langle \, , \, \rangle_p$ which to each $p \in M$ associates an inner product $\langle \, , \, \rangle_p$ on the tangent space at $p$, i.e. a positive-definite, symmetric, bilinear form. Further, the association $p \to \langle \, , \, \rangle_p$ is to be smooth in the following sense: if $\mathcal{X}$ and $\mathcal{Y}$ are $C^\infty$ vector fields on $M$, then $p \to \langle \mathfrak{X}(p), \mathfrak{Y}(p) \rangle_p$ is to be a $C^\infty$ function on $M$. One calls $M$ a *Riemannian manifold* and $\langle \; , \; \rangle_p$ a *Riemannian metric*. Note that given a smooth manifold $M$, it is not hard to prove the existence of a Riemannian metric on $M$: one simply embeds $M$ smoothly into some Euclidean space and then pulls back the standard inner product on $\mathbb{R}^n$.

One then applies a classical result of Levi-Civita establishing the existence of an *affine connection*. If we denote the algebra of $C^\infty$ vector fields on $M$ by $\mathcal{V}(M)$, an *affine connection* is a mapping $\nabla : \mathcal{V}(M) \times \mathcal{V}(M) \to \mathcal{V}(M)$, denoted by $(\mathfrak{X}, \mathfrak{Y}) \to \nabla_\mathfrak{X} \mathfrak{Y}$ which satisfies the following properties:

1. $\nabla_{f\mathfrak{X} + g\mathfrak{Y}} \mathfrak{Z} = f \nabla_\mathfrak{X} \mathfrak{Z} + g \nabla_\mathfrak{Y} \mathfrak{Z}$

2. $\nabla_\mathfrak{X} (\mathfrak{Y} + \mathfrak{Z}) = \nabla_\mathfrak{X} \mathfrak{Y} + \nabla_\mathfrak{X} \mathfrak{Z}$

3. $\nabla_\mathfrak{X} (f\mathfrak{Y}) = f \nabla_\mathfrak{X} \mathfrak{Y} + \mathfrak{X}(f) \mathfrak{Y}$

for any $\mathfrak{X}, \mathfrak{Y}, \mathfrak{Z} \in \mathcal{V}(M)$, and any $f, g \in C^\infty(M)$, see [17] or [5].

The presence of an inner product allows one to define the length of a $C^1$ curve $\gamma : [a, b] \to M$ by the formula $l(\gamma) = \int_a^b \sqrt{\langle \gamma'(t), \gamma'(t) \rangle_{\gamma(t)}} \, dt$ and subsequently define a metric $d : M \times M \to \mathbb{R}_{\geq 0}$ on $M$ via $d(x, y) = \inf l(\gamma)$, where the infimum is taken over all $C^1$ paths joining $x$ to $y$.

On the other hand, the presence of the Levi-Civita connection allows one to define the *curvature tensor* on $M$. The curvature tensor $R$ is a map which to each pair $(\mathfrak{X}, \mathfrak{Y}) \in \mathcal{V}(M) \times \mathcal{V}(M)$ associates a self-map $\mathcal{V}(M) \to \mathcal{V}(M)$ given by $R(\mathfrak{X}, \mathfrak{Y})\mathfrak{Z} = \nabla_\mathfrak{Y} \nabla_\mathfrak{X} \mathfrak{Z} - \nabla_\mathfrak{X} \nabla_\mathfrak{Y} \mathfrak{Z} + \nabla_{[\mathfrak{X}, \mathfrak{Y}]} \mathfrak{Z}$, for $\mathfrak{Z} \in \mathcal{V}(M)$. It turns out that the quantity $R(\mathfrak{X}, \mathfrak{Y})\mathfrak{Z}$ at $p \in M$ depends only on the values of $\mathfrak{X}$, $\mathfrak{Y}$, and $\mathfrak{Z}$ at $p$. Therefore, for any $x, y, z, t \in T_p M$, we can unambiguously define the number $(x, y, z, t)_p = \langle R(\mathfrak{X}, \mathfrak{Y})\mathfrak{Z}, \mathfrak{T} \rangle_p$ for any vector fields $\mathfrak{X}, \mathfrak{Y}, \mathfrak{Z}, \mathfrak{T} \in \mathcal{V}(M)$ which evaluated at $p$ give $x, y, z$ and $t$ respectively. Now, given any two-dimensional subspace $\sigma \subseteq T_p M$, we define the *sectional curvature* of $\sigma$ at $p$ as the number $K(\sigma) = \dfrac{(x, y, x, y)_p}{\sqrt{\langle x, x \rangle_p \langle y, y \rangle_p - \langle x, y \rangle_p}}$, where $\{x, y\}$ is any basis of $\sigma$. Again, one checks that this definition does not depend on the choice of a basis.

We are now finally ready to state Cartan's theorem:

**Theorem 2.2.7.** *(E Cartan) Let $M$ be a Riemannian manifold such that $K(\sigma) \leq 0$, for every 2-plane $\sigma \subseteq T_pM$. Then, $(M, d)$ is a CAT(0) space.*

*Proof.* See [6] or [3]. □

**Corollary 2.2.8.** *The hyperbolic spaces $\mathbb{H}^n$ are CAT(0).*

### New Spaces from Old

Two very useful constructions in topology are taking products and gluing. It turns out that if these are carried out in the right way, they will preserve the property of being non-positively curved in the sense of Definition 2.2.6.

**Proposition 2.2.9.** *Let $(X_1, d_1)$ and $(X_2, d_2)$ be CAT(0) spaces. Then the space $X = X_1 \times X_2$ with the metric $d((x_1, x_2), (y_1, y_2)) = \sqrt{d_1(x_1, y_1)^2 + d_2(x_2, y_2)^2}$ is also a CAT(0) space.*

Before we can state the next proposition, we need to recall the following basic metric construction. Let $(X_i, d_i)_{i \in I}$ be a family of metric spaces and let $C_i \subseteq X_i$ be a closed subspace of $X_i$ with each $C_i$ isometric to a given metric space, say $C$. Let $\iota_i : C \to C_i$ denote the isometry between $C$ and $C_i$. Then, the *amalgamation* of $X_i$ along $C$ is defined as the quotient of the disjoint union $\coprod_{i \in I} X_i$ by the equivalence relation generated by $(\iota_i(c) \sim \iota_j(c), \forall i, j \in I, c \in C)$. The amalgamation can be given a metric $d$ defined by:

- $d(x, y) = d_i(x, y)$, if $x, y \in X_i$ for some $i \in I$

- $d(x, y) = \inf_{c \in C} \{d_i(x, \iota_i(c)) + d_j(\iota_j(c), y)\}$, if $x \in X_i$, $y \in X_j$ with $i \neq j$

such that we have the following:

**Lemma 2.2.10.** *If $I$ is finite and each $X_i$ is proper, then $X$ is proper. If each $X_i$ is a geodesic space and $C$ is proper, then $X$ is a geodesic space.*

*Proof.* See Lemma 5.24 in [3] □

Fortunately for us, this well-known construction from topology preserves the property of being non-positively curved:

**Proposition 2.2.11.** *Let $X_1$ and $X_2$ be CAT(0) spaces and let $C$ be a complete metric space. Suppose that we have isometric embeddings $\iota_1 : C \to X_1$ and $\iota_2 : C \to X_2$ such that $\iota_1(C)$ and $\iota_2(C)$ are convex subsets of $X_1$ and $X_2$ respectively. Then, $X_1 \coprod_C X_2$ is a CAT(0) space.*

**Piecewise Euclidean CAT(0) Polyhedral Complexes**

Non-positively curved piecewise Euclidean polyhedral complexes are an important class of CAT(0) spaces and provide the first example for us of non-positively curved spaces which do not come from Theorem 2.2.7. In this section, we follow closely the exposition in I.7 of [3].

We shall call the convex hull of a finite set of points in $\mathbb{E}^n$ a *polyhedral cell*, or simply a *polyhedron*. Let $P \subset \mathbb{E}^n$ be a polyhedron, the dimension of the smallest affine subspace of $\mathbb{E}^n$ which contains $P$ is called the *dimension* of $P$. The interior of $P$ viewed as a subset of this affine plane is called the *interior* of $P$. If $H \subset \mathbb{E}^n$ is a hyperplane in $\mathbb{E}^n$, i.e. an affine subspace of $\mathbb{E}^n$ of codimension 1, then $H$ separates $\mathbb{E}^n$ into two closed half-spaces. If $P$ is contained in one of these closed half-space and $P \cap H \neq \emptyset$, we shall call $F = P \cap H$ a *face* of $P$, and if further $F \neq P$, a *proper face* of $P$. The support of $x \in P$ is defined to be the unique face of $P$ which contains $x$ in its interior.

**Definition 2.2.7.** (Piecewise Euclidean Polyhedral Complex) Let $\{P_i\}_{i \in I}$ be a family of polyhedra not necessarily of the same dimension. Let $X = \coprod_{i \in I} P_i$ denote their disjoint union. Let $\sim$ be an equivalence relation on $X$ and let $K = X/\sim$. Let $p : X \to K$ be the natural projection map, and let $p_i : P_i \to K$ denote the inclusion of $P_i$ into $X$ followed by $p$. We call $K$ a *piecewise Euclidean polyhedral complex* if:

1. for every $i \in I$, $p_i$ restricts to an injection on the interior of $P_i$,

2. for all $i, j \in I$ and $x_1 \in P_i$, $x_2 \in P_j$, if $p_i(x_1) = p_j(x_2)$, then there is an isometry $\iota : supp(x_1) \to supp(x_2)$ such that $p_i(y) = p_j(\iota(y))$ for all $y \in supp(x_1)$.

Definition 2.2.7 simply says that a polyhedral complex is an object obtained by gluing convex Euclidean polyhedra via linear isometries between their faces. Condition (1) guarantees that points in the interior of a polyhedron not get identified, while condition (2) guarantees that two faces be either identified via a linear isometry between them or not at all. A polyhedral complex all of whose polyhedra are 2-dimensional is called a *polygonal complex*. The set of isometry classes of faces of the polyhedra $P_i$ is denoted $Shapes(K)$.

We endow the polyhedral complex $K$ with the following pseudometric: $d(x, y)$ is defined to be the infimum of the lengths of piecewise linear paths $\lambda : [a, b] \to K$, where the length of a segment passing through a polyhedron is computed using the usual Euclidean metric on that polyhedron. A relevant for us special case of a theorem Bridson proved in his Ph.D. thesis, guarantees that under fairly mild hypotheses a piecewise Euclidean polyhedral complex is indeed a metric space:

**Theorem 2.2.12.** *If a piecewise Euclidean polyhedral complex has only finitely many isometry types of polyhedra, then it is a complete geodesic metric space.*

*Proof.* See [2] or [3]. □

For $x \in K$, we call the union of the interiors of the polyhedra that contain $x$, the *open star* of $x$, and denote it by $st(x)$. Now, let us fix $x \in K$. Given $y_1, y_2 \in st(x) - \{x\}$, we shall say that the geodesic segments $[x, y_1]$ and $[x, y_2]$ define the same *direction* at $x$ if either $[x, y_1] \subseteq [x, y_2]$, or $[x, y_2] \subseteq [x, y_1]$. The *link* of $x$ in $K$ is the set of directions at $x$. We shall denote the link of $x$ in $K$ by $Lk(x, K)$. The measure of a Euclidean angle carried by each polyhedron give rise to a pseudometric on $Lk(x, K)$. Under the hypotheses of Theorem 2.2.12, this pseudometric becomes a metric on $Lk(x, K)$.

Now that we have "made" the polyhedral complex into a complete geodesic metric space, we can ask for conditions which will guarantee that this metric space is CAT(0). We have thus far carefully avoided the discussion of the more general notion of a CAT($\kappa$) space, in order to keep the statements of the theorems simple. For the next result, however, we shall need the definition of a CAT(1) space. With the same setup as in Definition 2.2.6, we define a *CAT(1) metric space* $X$ to be a metric space in which every geodesic triangle $\Delta(a, b, c)$ having perimeter less than $2\pi$ satisfies the following *CAT(1) comparison inequality*: if $\overline{\Delta}(\overline{a}, \overline{b}, \overline{c})$ is a comparison triangle for $\Delta(a, b, c)$ in the unit sphere $\mathbb{S}^2$, then for any $p, q \in \Delta(a, b, c)$, we have $d(p, q) \leq d(\overline{p}, \overline{q})$. The metric on $\mathbb{S}^2$ is the angular metric described in Example 2.2.1. Compare this definition to Definition 2.2.6. Also, note the condition on the perimeter of $\Delta(a, b, c)$. We need this restriction on the size of the geodesic triangle $\Delta(a, b, c)$ in order to ensure that the corresponding comparison triangle $\overline{\Delta}(\overline{a}, \overline{b}, \overline{c}) \subset \mathbb{S}^2$ exists.

Now, we are able to state a necessary and nearly sufficient condition for a piecewise Euclidean polyhedral complex to be CAT(0):

**Definition 2.2.8.** A piecewise Euclidean polyhedral complex $K$ is said to satisfy *the link condition* if for every vertex $v \in K$, the link $Lk(v, K)$ is a CAT(1) space.

The link condition by itself is not sufficient as it only ensures non-positive curvature locally. The transitional step from local to global is made by the following well-known theorem:

**Theorem 2.2.13.** *A piecewise Euclidean polyhedral complex $K$ with $Shapes(K)$ finite is a CAT(0) space if and only if $K$ is simply connected and satisfies the link condition.*

*Proof.* See Theorem 5.2 in [3]. □

Having discussed CAT(0) spaces at some length as a generalization of the notion of non-positive sectional curvature in differential geometry, let us now consider the corresponding extension of the notion of negative sectional curvature. Again, we make vague allusions to CAT($\kappa$) spaces ($\kappa \leq 0$) defined in a manner completely analogous to Definition 2.2.6 with the only difference being that we take our comparison triangles from some 2-dimensional geometry of constant negative curvature, which is, of course, just $\mathbb{H}^2$ with a rescaled metric. While these do have their uses, it turns out that a more appropriate generalization was made by Gromov in [11]. Gromov's definition, like Alexandrov's definition of a CAT(0) space, makes use of the observation that in the presence of negative curvature the geodesic triangles become "thin". Unlike the definition of a CAT(0) space, however, Gromov's definition does not recourse to comparing geometric objects in the space, viz geodesic triangles, to objects which live in a different space. Thus, Gromov's definition of hyperbolicity is more robust and consequently many of the nice properties enjoyed by Gromov hyperbolic spaces and the groups which act on them are not shared by their CAT(0) analogues. This difference between CAT(0) groups and hyperbolic groups, however, is not seen as a flaw but rather as a complication which gives rise to interesting questions.

**Definition 2.2.9.** (Slim Triangles) Let $\delta > 0$. A geodesic triangle in a metric space is called $\delta$-slim if each of its sides is contained in a $\delta$-neighborhood of the union of the other two sides. A geodesic space $X$ is called $\delta$-*hyperbolic* if every geodesic triangle in $X$ is $\delta$-slim. We shall call a space simply *hyperbolic* if it is $\delta$-hyperbolic for some $\delta > 0$.

**Example 2.2.14.** The hyperbolic space $\mathbb{H}^n$ is $\delta$-hyperbolic.

By contrast with the notion of CAT(0) which has to do with non-positive curvature, none of the Euclidean spaces $\mathbb{E}^n$ are hyperbolic for $n > 1$. In fact, as the Flat Plane theorem below shows, the presence of isometrically embedded higher rank Euclidean spaces is one of the main obstruction to hyperbolicity.

**Theorem 2.2.15.** *(Flat Plane Theorem) Let $X$ be a proper CAT(0) space such that $Isom(X)$ acts on $X$ with a compact quotient. Then, $X$ is hyperbolic if and only if it does not contain a subspace isometric to $\mathbb{E}^2$.*

*Proof.* See Theorem 1.5, Chapter III.H in [3]. □

### 2.2.4   Isometries of CAT(0) Spaces

Let us introduce some standard terminology. An action of a group $\Gamma$ on a topological space $X$ is a homomorphism $\Phi : \Gamma \to Homeo(X)$, where $Homeo(X)$ is the group of self-homeomorphisms of $X$. One usually makes no mention of the homomorphism $\Phi$ and simply writes $\gamma \cdot x$ for $\Phi(\gamma)(x)$. For a subset $Y \subseteq X$, we shall write $\gamma \cdot Y$ for the image of $Y$ under $\gamma$, and $\Gamma \cdot Y$ for $\bigcup_{\gamma \in \Gamma} \gamma \cdot Y$. The action $\Phi$ is called:

- *faithful* if $Ker\ \Phi = \{1\}$,

- *free* if for every $x \in X$ and every $\gamma \in \Gamma - \{1\}$, we have $\gamma \cdot x \neq x$,

- *cocompact* if there exists a compact subset $C \subseteq X$ such that $X = \Gamma \cdot C$.

As we stated in Section 2.1, the objects of our study are discrete groups of isometries. By discreteness, we refer to the following fundamental notion of geometric group theory:

**Definition 2.2.10.** (Properly Discontinuous Action) Let $\Gamma$ be a group acting by isometries on a metric space $X$. The action is called *properly discontinuous* if for each $x \in X$ there exists a number $r > 0$ such that the set $\{\gamma \in \Gamma : \gamma \cdot B(x,r) \cap B(x,r) \neq \emptyset\}$ is finite.

The somewhat antiquated term "properly discontinuous" is meant to stand in contrast to the continuous action of a Lie group on a manifold. Henceforth, we shall instead use the term *proper* to mean properly discontinuous in the context of group actions.

There is a classical description of the isometries of a CAT(0) space, made by analogy with the classification of the isometries of $\mathbb{H}^2$, which we make use of in subsequent chapters, so we will recall the basics here. For the remainder of this section, $X$ will be a CAT(0) space. For $\gamma \in Isom(X)$, the *translation length* of $\gamma$ is defined to be the number $|\gamma| = \inf \{d(\gamma \cdot x, x) : x \in X\}$. If $|\gamma|$ is realized for some $x \in X$, we distinguish the cases:

- $|\gamma| = 0$, in other words $\gamma$ fixes a point $x \in X$. In this case $\gamma$ is called *elliptic*.

- $|\gamma| \neq 0$. In this case $\gamma$ is called *hyperbolic* or *axial*.

If $|\gamma|$ is not realized for any $x \in X$, then $\gamma$ is called *parabolic*. The set of points in $X$ where $d(\gamma \cdot x, x)$ attains its minimum is denoted $Min(\gamma)$. Whenever $Min(\gamma)$ is non-empty, $\gamma$ is called *semi-simple*. Thus, elliptic and hyperbolic isometries are semi-simple, while parabolic isometries are not. Since the distance function in a CAT(0) space is convex, see Proposition 2.2.4, the set $Min(\gamma)$ is a closed, convex, $\gamma$-invariant set. In fact, the action of a hyperbolic isometry $\gamma$ on $Min(\gamma)$ is particularly nice:

**Theorem 2.2.16.** *Let $X$ be a CAT(0) space.*

1. *An isometry $\gamma \in Isom(X)$ is hyperbolic if and only if there exists a geodesic line $\lambda : \mathbb{R} \to X$ on which $\gamma$ acts by a non-trivial translation, namely $\gamma \cdot \lambda(t) = \lambda(t+a)$. Further, the constant $a$ is actually equal to $|\gamma|$. The set $\lambda(\mathbb{R})$ is called an **axis** for $\gamma$.*

2. *If $X$ is a complete metric space and some power $\gamma^m$ is hyperbolic, then $\gamma$ is itself hyperbolic.*

3. *The axes for a hyperbolic isometry $\gamma$ are all parallel to each other and their union is $Min(\gamma)$, further $Min(\gamma)$ is isometric to a product $Y \times \mathbb{R}$ and the action of $\gamma$ on it is of the form $(y, t) \mapsto (y, t + |\gamma|)$, $y \in Y$, $t \in \mathbb{R}$.*

*Proof.* See Theorem 6.8 in [3]. □

Finally, we are able to make the definition of CAT(0) group:

**Definition 2.2.11.** (CAT(0) Group) Let $G$ be a group which acts properly and cocompactly on some CAT(0) metric space, then $G$ is called a *CAT(0) group*.

By analogy with Definition 2.2.11 above, we provisionally define a *δ-hyperbolic group* to be a group which acts properly discontinuously and by isometries on some δ-hyperbolic proper metric space with a finite diameter quotient. We postpone the official definition of a hyperbolic group until the next section where we discuss Cayley graphs.

### 2.2.5 Cayley Graphs and Quasi-isometries

**Definition 2.2.12.** (Quasi-isometry) Let $(X_1, d_1)$ and $(X_2, d_2)$ be two metric spaces. A map $f : X_1 \rightarrow X_2$ is a $(C, \epsilon)$-*quasi-isometric embedding* if $\frac{1}{C}d_1(x, y) - \epsilon \leq d_2(f(x), f(y)) \leq Cd_1(x, y) + \epsilon$. If in addition there exists a constant $A \geq 0$ such that every point in $X_2$ is in a $A$-neighborhood of the image of $f$, then $f$ is called a $(C, \epsilon)$-*quasi-isometry*.

Note that neither the definition of a quasi-isometric embedding nor the definition of a quasi-isometry requires it to be a continuous map. It turns out, however, that quasi-isometries are the right kind of maps to consider if one is interested in the coarse structure of a space. In fact, the notion of a quasi-isometry is indispensable to any geometric group theorist as it will be to us in the remainder of this work.

Our point of view thus far has been this: we started out with a metric space, then considered subgroups of its isometry group. This way our groups came naturally with an action on some topological space. We shall now consider the dual point of view: we start with a finitely generated group $G$ then proceed to construct an action of $G$ on some topological space. In fact, it is an interesting question to find the "nicest" possible action that a given finitely generated group admits.

**Definition 2.2.13.** (Cayley Graph) Let $G$ be a finitely generated group, and let $\mathcal{S}$ be a finite set of generators for $G$. The Cayley graph of $G$, relative to the generating set $\mathcal{S}$, which we shall denote by $\Gamma_{\mathcal{S}}(G)$ is defined as follows: We take the group $G$ itself to be the set of vertices, $\mathcal{V} = G$. Two vertices $g_1, g_2 \in \mathcal{V}$ are joined by an edge if $g_2 = g_1 s$, for some $s \in \mathcal{S}$. In other words, $\mathcal{E} = \{(g, s) : g \in G, s \in \mathcal{S}\}$.

When endowed with the *length metric*, the Cayley graph becomes a metric space: we first declare each edge in $\Gamma_{\mathcal{S}}(G)$ to have length equal to 1, then define $d_{\mathcal{S}}(x, y)$ to equal the infimum of the lengths of paths joining $x$ to $y$. Since the group $G$ is now embedded as the vertex set of the metric graph $(\Gamma_{\mathcal{S}}(G), d)$, $G$ itself becomes a metric space. The metric $d_{\mathcal{S}}$ on $G$ is called the *word metric* on $G$ relative to the generating set $\mathcal{S}$. The reason for this terminology is the following:

**Lemma 2.2.17.** *Let $G$ be a finitely generated group with a finite symmetric generating set $\mathcal{S}$ (i.e. $\mathcal{S} = \mathcal{S}^{-1}$), and let $d$ be the metric on $G$ defined above. Given a word on the generators $w = s_1 ... s_n$ in $G$, let $\|w\| = n$. Then, for $g_1, g_2 \in G$, $d_{\mathcal{S}}(g_1, g_2) = \inf \left\{ \|w\| : w = g_1^{-1} g_2 \in G \right\}$.*

*Proof.* Obvious. □

Definition 2.2.13 begs the question: how does the word metric depend on the generating set? The answer, which is simple enough, is given by:

**Lemma 2.2.18.** *Let $G$ be a finitely generated group and let $\mathcal{S}_1$ and $\mathcal{S}_2$ be two finite symmetric generating sets for $G$. Then the identity map $\iota : (G, d_{\mathcal{S}_1}) \to (G, d_{\mathcal{S}_2})$ is a quasi-isometry. The analogous statement for Cayley graphs is also true: $\iota : \Gamma_{\mathcal{S}_1}(G) \to \Gamma_{\mathcal{S}_2}(G)$ is a quasi-isometry, where $\iota$ is any map which extends the identity map on $G \subseteq \Gamma_{S_i}(G)$, for $i = 1, 2$.*

*Proof.* Exercise. □

Due to Lemma 2.2.18, one usually suppresses the mention of $\mathcal{S}$. Since the edges emanating from a vertex are labeled by elements in the generating set, we easily conclude that $\Gamma(G)$ is a proper metric space. Further, the group $G$ acts on $\Gamma(G)$ on the left via $\gamma \cdot (g, s) = (\gamma g, s)$. It is easily seen that this action is by isometries. We are now ready to state the definition of a hyperbolic group:

**Definition 2.2.14.** (Hyperbolic Group) A finitely generated group $G$ is called a hyperbolic group if its Cayley graph is a hyperbolic metric space.

We note that Definition 2.2.14 makes sense in view of:

**Theorem 2.2.19.** *The property of being hyperbolic is a quasi-isometry invariant: If $X$ is a $\delta$-hyperbolic metric space, then any space to which $X$ is quasi-isometric is also $\delta$-hyperbolic (with a possibly different value of $\delta$).*

*Proof.* See Theorem 1.9, Chapter III.H in [3]. □

## 2.3 Ends of Pairs of Groups and Sageev's Cubing

### 2.3.1 Overview of Bass-Serre Theory

A central idea in combinatorial group theory is to investigate the algebraic structure of groups which act on trees. First, let us define some basic notions. Given a group $G$ and a subset $S \subseteq G$, the *normal closure* $\langle\langle S \rangle\rangle$ of $S$ in $G$ is defined to be the smallest normal subgroup of $G$ which contains $S$.

**Definition 2.3.1.** (Free Product with Amalgamation) Let $G_1$, $G_2$, and $H$ be groups, and let $\iota_1 : H \to G_1$, $\iota_2 : H \to G_2$ be monomorphisms. The *amalgamated product* of $G_1$ and $G_2$ along $H$ is the group $G_1 * G_2 / \langle\langle \iota_1(h)\iota_2(h^{-1}) \rangle\rangle_{h \in H}$ and we shall denote it by $G_1 *_H G_2$.

One may be tempted to define the amalgamated product as a pushout of a certain diagram in the category of groups and then show its existence. However, this approach does not give any idea as to the group structure and since we have no need for categories, we much prefer this explicit description instead.

**Definition 2.3.2.** (HNN Extension) Let $G$ and $H$ be groups and let $\iota_1, \iota_2 : H \to G$ be two monomorphisms. If $S$ is any generating set for $G$ and $t \notin S$, the group defined by the presentation $\langle S, t : t\iota_1(h)t^{-1} = \iota_2(h) \rangle$ is called the *HNN extension* of $G$ over $H$ relative to $\iota_1$ and $\iota_2$.

Amalgamated products and HNN extensions arise naturally as the fundamental groups of gluings of topological spaces.

**Example 2.3.1.** Let $(X_1, x_1)$, $(X_2, x_2)$, and $(Y, y_0)$ be three based topological spaces. Let $f_i : (Y, y_0) \to (X_i, x_i)$, for $i = 1, 2$, be continuous injections which induce monomorphisms of fundamental groups $f_{i*} : \pi_1(Y, y_0) \to \pi_1(X_i, x_i)$, for $i = 1, 2$. Then, $\pi_1(X_1 \coprod_Y X_2, [x_0]) \cong \pi_1(X_1, x_1) *_{\pi_1(Y, y_0)} \pi_1(X_2, x_2)$ by the Seifert-van Kampen theorem.

If instead of gluing two distinct spaces along a common subspace, one glues two subspaces $Y_i$ of the topological space $(X, x_0)$ via a homeomorphism between them, one ends up with a fundamental group which is isomorphic to an HNN extension of the fundamental group of $(X, x_0)$. We refer the interested reader to [31] where the details of these computations have been worked out.

Whenever a group $G$ is isomorphic to either a non-trivial free product with amalgamation or an HNN extension over $H$, we shall say that $G$ *splits* over $H$. As it turns out, there is a close connection between group splittings and actions on trees as evidenced by the "fundamental theorem" of Bass-Serre theory:

**Theorem 2.3.2.** *Let $G$ be a finitely generated group. Then, $G$ splits as a non-trivial free product with amalgamation, or as an HNN extension, if and only if $G$ acts on a* *tree without a global fixed point.*

*Proof.* See [31]. □

In fact, one can explicitly construct the tree and the group action on it, so this is not purely an existence result. Conversely, if such an action exists, then one constructs a decomposition of the group $G$ into free products with amalgamation and HNN extensions over subgroups which appear as vertex and edge stabilizers of the tree action.

Perhaps even more remarkable is the classification theorem of finitely generated groups according to their number of ends, see Theorem 2.3.3. The reader might be familiar with the notion of *ends* of a metric space. Let $(X, d)$ be a metric space, and let $Y \subseteq X$ be a subspace. Let us denote the number of connected components of $Y$ which are unbounded by $c(Y)$. Then, the number of ends $e(X)$ of $X$ is defined as $e(X) = \sup \{c(X - K) : K \subseteq X\}$, where the supremum is take over all compact subsets of $X$. Informally, $e(X)$ is the number of ways to approach infinity in $X$. We can now unceremoniously define the number of ends of a group $G$ to be the number of ends of its Cayley graph $\Gamma(G)$. However, instead of doing this and then having to check that the definition is indeed independent of the generating set, we shall outline the more constructive definition found in [31].

For any set $S$, the *power set* of $S$ denoted by $\mathcal{P}(S)$ is the set of subsets of $S$. It is an exercise in an introductory class on set theory or real analysis to make $\mathcal{P}(S)$ into an abelian group under the operation of taking symmetric difference: $S_1 + S_2 = S_1 \Delta S_2$. The identity in this group is the empty set, and every element is its own inverse. Now, if $G$ is a group, we see that the subset $\mathcal{F}(G)$ of $\mathcal{P}(G)$ consisting of finite subsets of $G$ forms a subgroup of $\mathcal{P}(G)$. So does the subset $\mathcal{Q}(G) = \{A \in \mathcal{P}(G) : card\,(A + Ag) < \infty \; \forall g \in G\}$. The elements of $\mathcal{Q}(G)$ are called *almost invariant* sets. Now, note that since every element in the abelian group $\mathcal{P}(G)$ has order equal to 2, $\mathcal{P}(G)$ is a vector space over $\mathbb{Z}_2$. The same is true for $\mathcal{F}(G)$ and

$\mathcal{Q}(G)$, respectively.

**Definition 2.3.3.** Let $G$ be a finitely generated group. The *number of ends* of $G$ is the number $e(G) = \dim_{\mathbb{Z}_2}(\mathcal{Q}(G)/\mathcal{F}(G))$.

It turns out that unlike the number of ends of a topological space which can equal really anything, the values of $e(G)$ are severely restricted:

**Theorem 2.3.3.** *Let $G$ be a finitely generated group. Then, $e(G) = 0, 1, 2$ or $\infty$ and:*

1. *$e(G) = 0$ if and only if $G$ is finite,*

2. *$e(G) = 2$ if and only if $G$ contains $\mathbb{Z}$ with finite index,*

3. *if $e(G) = \infty$, then $G$ splits over a finite subgroup.*

*Proof.* See [31]. □

The proof of part (1) of Theorem 2.3.3 is not too difficult and follows from the basic properties of ends outlined in [31]. Part (2) is due to Hopf along with the observation that $e(G)$ can only equal one of $0, 1, 2, \infty$, while part (3) is a celebrated result due to John Stallings. In fact, parts (2) and (3) of Theorem 2.3.3 have the following converse of sorts:

**Theorem 2.3.4.** *If $G$ splits over a finite subgroup, then $e(G) \geq 2$.*

*Proof.* See [31]. □

As a basic property of the number of ends of a group, we have:

**Theorem 2.3.5.** *If $G$ acts freely on a connected CW-complex $X$ such that the quotient $K$ is a finite CW-complex, in other words $G$ is the group of automorphisms of the cover $X \to K$, then $e(G) = e(X)$.*

*Proof.* See Theorem 5.4 in [31]. □

As a special case of Theorem 2.3.5, we have the following result which justifies our initial informal definition of $e(G)$:

**Proposition 2.3.6.** *Let $G$ be a finitely generated group, and let $S$ be any finite set of generators for $G$. Then, the number of ends of $G$ is the same as the number of ends of the Cayley graph $\Gamma_S(G)$.*

*Proof.* See Proposition 5.2 in [31]. $\qquad\square$

### 2.3.2 Ends of Group Pairs

In this section we generalize the notion of ends of a group. Given a subgroup $H \leq G$, we shall define the number of ends of the pair $(G, H)$ in such a way that when $H$ is the trivial group $e(G, H) = e(G)$. We shall again follow the exposition in [31]. In order to better understand the definition of $e(G, H)$, let us recast the definition of $e(G)$ in a new light. Again, we start with a finitely generated group $G$, and define the $\mathbb{Z}_2$-vector spaces $\mathcal{P}(G)$ and $\mathcal{F}(G)$ as before. Let $\mathcal{E}(G) = \mathcal{P}(G)/\mathcal{F}(G)$. The right action of $G$ on itself by right-multiplication induces an action of $G$ on $\mathcal{E}(G)$. If we let $\mathcal{E}(G)^G$ denote the set of fixed points of this action, then with the notation already established in Section 2.3.1, we have $\mathcal{E}(G)^G = \mathcal{Q}(G)/\mathcal{F}(G)$. Thus, our definition of $e(G)$ becomes $e(G) = \dim_{\mathbb{Z}_2} \mathcal{E}(G)^G$. Now, let $H$ be a subgroup of $G$, and let $H \backslash G$ be the set of right cosets of $H$. We define the *number of ends of the pair* $(G, H)$ to be the number $e(G, H) = \dim_{\mathbb{Z}_2} \mathcal{E}(H \backslash G)^G$. The elements of $\mathcal{E}(H \backslash G)^G$ are called *H-almost invariant sets*, and the nonzero *H-almost invariant sets* that are not in the same equivalence class as $G$ are called *proper*. The analogue of Theorem 2.3.5 is the following:

**Lemma 2.3.7.** *Let $\tilde{X} \to X$ be a regular cover of the finite CW-complex $X$ with automorphism group $G$. If $H$ is a subgroup of $G$, then $e(G, H) = e(H \backslash X)$.*

*Proof.* See Lemma 8.1 in [31]. $\qquad\square$

One even has the following counterpart to Theorem 2.3.4:

**Lemma 2.3.8.** *If $G$ splits over a subgroup $H \leq G$, then $e(G, H) \geq 2$.*

*Proof.* See Lemma 8.3 in [31] □

Unfortunately, this is as far as the analogy can be pushed. Namely, an example due to Scott exhibits an unsplittable group $G$ and a subgroup $H \leq G$ such that $e(G, H) \geq 2$. Again, see Chapter 8 of [31]. Thus, in light of Theorem 2.3.2, there is no action of $G$ on a tree without a fixed point. Now, in the absence of "interesting" actions on trees, one may lose hope of ever having a nice geometrical portrait of $G$ in the style of Bass-Serre. This is where Sageev's celebrated construction comes to the rescue.

### 2.3.3 Sageev's Cubing

**Overview**

The remark at the end of Section 2.3.2 reveals that the case when the group pair $(G, H)$ has at least two ends merits further investigation. We shall call such a group pair *multi-ended* following [26]. Equivalently, we say that $H$ is a *codimension-1* subgroup of $G$. The terminology "codimension-1" comes from topology, where such a subgroup is typically obtained as the image of the fundamental group of a manifold under the map induced by a $\pi_1$-injective immersion into a manifold of one higher dimension. If a group $G$ has a subgroup $H$ such that the pair $(G, H)$ is multi-ended, $G$ will be called *semi-splittable*. In this section, we present a summary of Sageev's cubing construction and we closely follow the exposition in [26]. A *cube complex $X$* is a polyhedral complex in which each polyhedron is isometric to the unit cube in $\mathbb{E}^n$, for some $n$. If a cube complex given the path metric described following Definition 2.2.7 is a CAT(0) space, $X$ is called a *cubing*. Cubings are a natural generalization of trees in that a tree is simply a 1-dimensional cubing. The statement of Sageev's result is a natural generalization of Theorem 2.3.2:

**Theorem 2.3.9.** *Let $G$ be a finitely generated group. Then, $G$ is semi-splittable if and only if $G$ acts essentially on a cubing.*

*Proof.* See Theorem 3.1 in [26]. $\square$

A few words are in order regarding the word "essentially" in the theorem above. The existence of an action of $G$ is not itself surprising as we can always consider the trivial action of $G$ on any space. Thus, we need to set a standard which an action will need to satisfy before it can be considered "interesting". Theorem 2.3.2 asserts the existence of an interesting action on a tree $T$ in the sense that no point of $T$ is fixed by all $g \in G$. It turns out that an action of group on a tree without a global fixed point implies the existence of $g \in G$ which as an automorphism of the tree has no fixed points, see [26]. In order to define properly what an essential action is, we need to talk about the geometry of a cubing.

Let $X$ be a cubing and let $\mathcal{E}$ be the set of oriented edges of $X$. We shall define two edges $e, f \in \mathcal{E}$ to be equivalent if there exists a finite sequence of edges $f = e_0, ..., e_n = e$, such that $e_i$ and $e_{i+1}$ are opposite sides of some 2-cube in $X$. This obviously defines an equivalence relation on $\mathcal{E}$ and we shall call an equivalence class in $\mathcal{E}$ under this relation a *combinatorial hyperplane*. Let us explain the reason for this terminology. Recall that if $\sigma$ is a $n$-cube and $e$ an edge of $\sigma$, then the cube dual to $e$ in $\sigma$ is the intersection of $\sigma$ with a hyperplane in $\mathbb{E}^n$ which is orthogonal to $e$ and bisects it. Note that dual cubes live combinatorially in the first barycentric subdivision of a cube complex. Now, let $J$ be the subcomplex of the first barycentric subdivision of $X$ which consists of cubes dual to the edges in $H$. This subcomplex is called the *geometric realization* of the combinatorial hyperplane $H$ or simply a *geometric hyperplane*. A geometric hyperplane $J$ is said to *self-intersect* if some cube of $X$ contains two cells belonging to $J$. It is a fundamental fact that geometric hyperplanes in cubings behave much like hyperplanes in $\mathbb{E}^n$, namely:

**Theorem 2.3.10.** *Suppose that $X$ is a cubing and $J$ a geometric hyperplane in $X$. Then, $J$ does not self-intersect and $X - J$ has exactly two components.*

*Proof.* See [7]. □

The *oriented stabilizer* of $J$, which we shall denote by $stab(J)$, is the subgroup of $G$ which consists of elements of $G$ which preserve the partition of $X$ by $J$ setwise. In other words, if $Y$ and $Y^c$ denote the two components of $X - J$, then $stab(J) = \{g \in G : g \cdot Y = Y\}$. We shall say that the action of $G$ is *essential with respect to $J$* if there exists a vertex $v$ in $X$ such that both $V = \{g \in G : g \cdot v \in Y\} \subseteq G$ and $V^c = G - V$ contain infinitely many right cosets of $stab(J)$. Finally, we call an action *essential* if it is essential with respect to some hyperplane.

**The Cubing Construction**

We now describe in more detail how one constructs the cubing promised in Theorem 2.3.9. Let $G$ be a finitely generated group and let $H$ be a codimension-1 subgroup. This means that we have a proper $H$-almost invariant set $A \subset G$. By modifying $A$ if necessary, we may assume that $A$ is an *invariant $H$-set*, meaning that $A$ is $H$-almost invariant and $h \cdot A = A$, for all $h \in H$, see Section 2.2 of [26]. Now, we consider the set $\Sigma = \{gA : g \in G\} \cup \{gA^c : g \in G\}$ partially ordered under inclusion. First, we construct the vertices of the cubing. The vertices of $X$ will be the subsets $V$ of $\Sigma$ which satisfy the following two properties:

- Given $A \in \Sigma$, then either $A \in V$, or $A^c \in V$, but not both.

- If $A \in V$ and $A \subseteq B$, then $B \in V$.

Let us denote this set of vertices by $\tilde{\mathcal{V}}$. Upon seeing this construction for the first time, one cannot help but notice that the elements of the vertex set look very much like ultrafilters on $\Sigma$. We recall here that an *ultrafilter* $\mathcal{U}$ on a set $S$ is a subset of $\mathcal{P}(S)$ such that:

1. $\emptyset \notin \mathcal{U}$,

2. if $A, B \in \mathcal{P}(S)$, and $A \in \mathcal{U}$, then $A \subseteq B$ implies that $B \in \mathcal{U}$,

3. for any $A \in \mathcal{P}(S)$, precisely one of $A$ or $A^c$ belongs to $\mathcal{U}$,

4. if $A, B \in \mathcal{U}$, then $A \cap B \in \mathcal{U}$.

The only difference between an ultrafilter on $\Sigma$ and a vertex is that a vertex need not satisfy (4). Of the set of all ultrafilters on $S$, one distinguishes those which are generated by a single element $s$ of $S$. More precisely, consider the ultrafilters $\mathcal{U}_s = \{A \subseteq S : s \in A\}$. These are called *principal ultrafilters*, perhaps by analogy with the notion of principal ideals in algebra, and they correspond to the following vertices:

**Example 2.3.11.** For $g \in G$, the subset $V_g = \{A \in \Sigma : g \in A\}$ is easily verified to be a vertex. In fact, the vertices of this form are precisely those which are also principal ultrafilters on $\Sigma$. They are called *basic vertices*.

Now, we join two vertices $V_1$ and $V_2$ of $X$ by an edge if, as subsets of $\Sigma$, they differ by exactly one element, i.e. $V_2 = (V_1 - \{A\}) \cup A^c$, for some $A \in V_1$. This set of edges we shall denote by $\tilde{\mathcal{E}}$. We now have constructed a 1-dimensional complex, namely the graph $(\tilde{\mathcal{V}}, \tilde{\mathcal{E}})$ which is too large in the sense that it not only contains a "copy" of $G$ comprised of the vertices $V_g$ described in Example 2.3.11 but a lot of other stuff as well. We shall therefore trim down this graph as follows: let $\mathcal{V}$ be the set of vertices $V$ in $\tilde{\mathcal{V}}$ for which there exists an edge path from $V$ to $V_g$ for some $g \in G$. We shall denote the set of vertices and edges of the resulting graph by $X^{(1)} = (\mathcal{V}, \mathcal{E})$. Once the 1-skeleton has been constructed, one proves that it is connected, see Theorem 3.3 in [26]. Then, we begin to attach higher dimensional cubes inductively: we attach an $(n + 1)$-cube if we see its boundary in $X^{(n)}$.

Next, we record some observations all found in [26] leading up to Theorem 2.3.14 below, which will be interesting to us in their own, thereby reversing their logical order.

**Lemma 2.3.12.** *Suppose $V$ is a vertex and $A \in V$. Then $W = (V - \{A\}) \cup \{A^c\}$ is also a vertex if and only if $A$ is not properly contained in a subset $B \in V$, in other words if and only if $A$ is **minimal** in $V$ with respect to inclusion.*

*Proof.* See Lemma 3.2 in [26]. $\square$

In light of Lemma 2.3.12, if $A$ is minimal in $V$, we shall denote the set $(V - A) \cup \{A^c\}$ by $(V; A)$ thus borrowing notation from [26]. We define recursively, $(V; A_1, ..., A_n) = \{(V; A_1..., A_{n-1}) - A_n\} \cup \{A_n^c\}$, provided that $A_n$ is minimal in $(V; A_1, ..., A_{n-1})$. The ordering of $A_1, ..., A_n$ is immaterial. This is obvious, as we can think of the $A_i$ as certain marked elements of $V$ and note that different arrangements of the $A_i$ give rise to the same marking of $V$.

**Lemma 2.3.13.** *Suppose $\sigma^n$ is an $n$-cube in $X$ having $V$ as a vertex and suppose that the vertices of $\sigma^n$ connected to $V$ by an edge are $(V; A_1)$, $(V; A_2),...,(V; A_n)$. Then, the vertex opposite to $V$ is $(V; A_1, ..., A_n)$.*

*Proof.* See Lemma 3.5 in [26]. $\square$

With the notation in Lemma 2.3.13, if $V$ is a vertex of the $n$-cube $\sigma^n$ and $V' = (V; A_1, ..., A_n)$ is its opposite vertex, then any other vertex of $\sigma^n$ is of the form $W = (V; A_{i_1}, ..., A_{i_k})$, for some subset $\{A_{i_1}, ..., A_{i_k}\} \subseteq \{A_1, ..., A_n\}$. The result, which properly belongs to polyhedral combinatorics, is a by-product of the proof of Lemma 2.3.13.

Finally, we state the central result of Sageev's paper [26]:

**Theorem 2.3.14.** *The set $X$ is a cubing.*

*Proof.* See Theorem 3.7 in [26]. □

The statement of Theorem 2.3.14 would be a fitting conclusion of this chapter were it not for the following fact: Sageev's cubing need not be finite dimensional. In this case, the action of $G$ on $X$ has no hope of being cocompact. To a geometric group theorist, this is a rather upsetting conclusion. Further, the cubing, by its very nature, depends on the choice of an $H$-almost invariant set $A$, and we shall henceforth use $X_A$ to emphasize this.

# Chapter 3

# On the Diagonal Action of $\mathbb{Z}$ on $\mathbb{Z}^n$

Most of the work of this chapter is contained in [28], however, an important part of it will appear for the first time in this dissertation. In this chapter we study the action of $\mathbb{Z}$ on $\mathbb{Z}^n$ from a dynamical perspective. The motivation for this study comes from the notion of bounded packing introduced by Hruska and Wise in [14]. Our work in this chapter, however, has only partial overlap with the original idea in [14] and has some interesting consequences in terms of growth of the set of left cosets which we study in Section 3.7.2.

## 3.1 The Origin of Bounded Packing

### 3.1.1 Disc Packing in the Plane

The term *bounded packing* was first used by Hruska and Wise in [14]. Very informally, if $G$ is a finitely generated group and $H$ a subgroup, $H$ is said to have bounded packing if the cosets of $H$ are not densely packed together. Before we make the notion precise let us underline the geometrical nature of the idea. We become familiar with packing problems from early on in life. For example, take the plane $\mathbb{E}^2$ and consider the collection of closed discs $\mathcal{S} = \left\{ B(a, \frac{1}{2}) : a \in \mathbb{Z}^2 \subseteq \mathbb{E}^2 \right\}$. Very informally, we shall say that the discs in $\mathcal{S}$ are packed because no two of them overlap except along their boundary. Note that every point in $\mathbb{E}^2$ is within finite distance of some disc in $\mathcal{S}$. This way of packing $\mathbb{E}^2$ is very nice in the sense that the number of discs in $\mathcal{S}$ which intersect any given ball of radius $r > 0$ is finite, in other words $\mathcal{S}$ has the property of bounded packing in $\mathbb{E}^2$. It is not too hard to enlarge $\mathcal{S}$ so as to obtain a new collection of discs packing $\mathbb{E}^2$ which does not have the local finiteness property - we only need to keep inserting discs in the spaces enclosed by discs in $\mathcal{S}$ such that each newly inserted disc has a boundary circle tangent to the boundary circles of all the existing discs. We keep "filling in the gaps" in this manner, a task which we easily convince ourselves can be accomplished, to build a new collection of discs $\mathcal{S}'$ which does not have bounded packing.

In the context of group theory however, bounded packing has become intimately connected to Sageev's cubing via Lemma 3.1.5 below which is found in Sageev's paper [27], and a special case of which is attributed by him to Hopf.

### 3.1.2 More on Cubings

Let $H$ be a codimension-1 subgroup of the finitely generated group $G$. This means that $e(G, H) \geq 2$, which by definition means that we have a proper $H$-almost invariant set $A$, see Section 2.3.2. Recall that in order to construct the Sageev cubing in Section 2.3.3, we worked with the set $\Sigma$ which consisted of left translates of $A$ and $A^c$. We shall say that two elements $A_1, A_2 \in \Sigma$ are *nested* if one of the following inclusions holds: $A \subset gA$, $A \subset gA^c$, $A^c \subset gA$, $A^c \subset gA^c$, following [27]. Further, we define $w(\Sigma)$, the *width* of $\Sigma$, to be the size of the largest collection of pairwise non-nested elements of $\Sigma$. We are interested in the width of $\Sigma$ because of the following result which in the chain of logical inference precedes Theorem 2.3.14:

**Lemma 3.1.1.** *Suppose $V$ is a vertex and $S = \{A_1, ..., A_n\} \subset V$. Then $(V, S)$ spans a cube $\sigma^n$ if and only if $A_i \neq A_j$, each $A_i$ is minimal in $V$, and no $A_i^c$ is contained in an $A_j$, for all $i, j \in \{1, ..., n\}$.*

*Proof.* For the proof, which essentially follows from Lemmas 2.3.12 and 2.3.13, see [26] Lemma 3.6. □

As an immediate consequence of Lemma 3.1.1 we have:

**Corollary 3.1.2.** *The dimension of $X_A$ is bounded above by the width of $\Sigma$.*

*Proof.* Obvious. □

We now seriously consider the question of the dimension of the Sageev cubing. In order to do that, we need the following definition:

**Definition 3.1.1.** Let $G$ be a finitely generated hyperbolic group with a generating set $\mathcal{A}$. The subgroup $H \subseteq$ is a *quasiconvex* subgroup of $G$ if $H$ is a quasiconvex subset of $\Gamma_{\mathcal{A}}(G)$.

Standard arguments show that this definition is independent of the generating set $\mathcal{A}$.

Even though we only define quasiconvexity here in order to state Theorem 3.1.4 below, we emphasize that this notion is central in geometric group theory, interesting in itself due to the following fact:

**Lemma 3.1.3.** *Let $G$ be a finitely generated group with a generating set $\mathcal{A}$ and let $H$ be a quasiconvex subgroup of $G$. Then, $H$ is finitely generated and the inclusion map $\iota : H \hookrightarrow G$ is a quasi-isometric embedding.*

Quasiconvexity will also be important for us in Chapter 4.

In general, it is possible, in fact even likely that the word length of an element $h \in H$ may decrease when written in the generating set of $G$, simply because we have a larger generating set to work with. In fact, the word length may shrink quite a bit. Analyzing this phenomenon, one is naturally led to the notion of *subgroup distortion*, which is implicit in the work in Section 3.7.2. What Lemma 3.1.3 says is that quasiconvex subgroups are undistorted in $G$.

We are now able to state the interesting result of Sageev's, which explains the relationship between quasiconvexity and his cubing construction:

**Theorem 3.1.4.** *Suppose that $G$ is a hyperbolic group and $H$ is a codimension-1 quasiconvex subgroup. Then, there exists an $H$-almost invariant set $A$ such that the cubing $X_A$ is finite dimensional and the action of $G$ on $X_A$ is cocompact.*

*Proof.* See Theorem 3.1 in [27]. □

It is in the proof of this theorem that Sageev implicitly makes use of the notion of bounded packing of the cosets of a certain subgroup of $G$. In order to see how the geometric notion of packing enters the combinatorial construction of Sageev, we need to turn our attention to the idea of ends once again.

Recall that according to the definition we gave in Section 2.3.2, $H$ is codimension-1, or equivalently the pair $(G, H)$ is multi-ended, if $G$ contains a proper $H$-almost invariant set $A$. Unfortunately, even though we know that such a set exists, we have no idea how to find it. At the same time, understanding how $A$ meets its translates $gA$ and the translates of its complement $gA^c$ in $G$ is essential to understanding the cubing construction; we therefore briefly describe how to explicitly construct such a proper $H$-almost invariant set in $G$, following [27].

Consider the Cayley graph $\Gamma(G)$ of $G$, and consider also the action of $H$ on it given by restricting the left action of $G$ on $\Gamma(G)$ described in Section 2.2.5. In view of Lemma 2.3.7 $e(G, H) = e(H\backslash\Gamma(G))$, hence some bounded neighborhood of $H$ in $H\backslash\Gamma(G)$ separates $H\backslash\Gamma(G)$ such that at least two components are unbounded. Therefore, after lifting this neighborhood to $\Gamma(G)$ we are able to find a bounded neighborhood $N_\nu(H)$ of $H$ in $\Gamma$ which separates $\Gamma$ in a way that at least two of the components of $\Gamma(G) - N_\nu(H)$ contain vertices arbitrarily far away from $H$. Now, the vertices of any one of these components can be taken to be our proper $H$-almost invariant set $A$, see Theorem 2.3 in [26].

In the setting of Theorem 3.1.4, $H$ is finitely generated which is a consequence of its being quasiconvex. This allows us to replace the original bounded neighborhood of $H$ with a connected one, which we still denote by $N_\nu(H)$, for a possibly larger value of $\nu$. As Lemma 3.1.1 above shows, the existence of an $n$-cube in $X_A$ can be detected by the presence of a subset of $\Sigma$ which consists of nonnested elements. Now, Lemma 3.1.5 below shows that the existence of large subsets of $\Sigma$ consisting of pairwise non-nested elements can be prevented by restricting the manner in which the left translates of $N_\nu(H)$ intersect one another.

**Lemma 3.1.5.** *If $gN_\nu(H) \cap N_\nu(H) = \emptyset$ for some $g \in G$, then $gA$ and $A$ are nested.*

*Proof.* See Lemma 3.2 in [27]. □

Given a metric space $(X, d)$, there is a well known distance function on $\mathcal{P}(X)$ defined as follows: $d(Y_1, Y_2) = \inf \{d(y_1, y_2) : y_1 \in Y_1, y_2 \in Y_2\}$ for $Y_1, Y_2 \subseteq X$. In light of the previous Lemma, if $gN_\nu(H) \cap N_\nu(H) \neq \emptyset$, then $d(gH, H) < 2\nu$, where $d$ is the distance function on $H \backslash G$ inherited from the word metric on $G$ defined in 2.2.5. We now come to the result in [27] which motivates Definition 3.2.1 below:

**Corollary 3.1.6.** *If every collection of left cosets of $H$ in $G$ containing at least $n$ elements has a pair of cosets $aH, bH$ with $d(aH, bH) > 2\nu$, then the width of $\Sigma$, and therefore $\dim(X_A)$ is bounded above by $n$.*

*Proof.* By Lemma 3.1.1, in order to have an $m$-dimensional cube in $X_A$, we must have a subset of $\Sigma$ consisting of $m$ pairwise nonnested elements. A subset of $\Sigma$ having this property corresponds to a collection of $m$ left translates of $N_\nu$ which pairwise intersect nontrivially; this is Lemma 3.1.5. Hence, we find a collection of $m$ left cosets of $H$ which are pairwise within distance $2\nu$ of each other. However, for $m > n$ such a collection of left cosets of $H$ cannot exist. □

## 3.2 Bounded Packing and Coset Growth - Definition and Basic Properties

If $G$ is a topological group, a metric $d$ on $G$ is called *left-invariant* if for all $g \in G$, the map $a \mapsto ga$ is an isometry of $G$ to itself. If $G$ is a finitely generated group, the word metric on $G$ corresponding to a finite generating set is clearly a left-invariant metric. Throughout this section, $G$ will be a finitely generated group, $d$ will be the word metric unless otherwise stated, and $H$ will be a subgroup.

**Definition 3.2.1.** (Bounded Packing) We say that $H$ has *bounded packing* in $G$, if given any $r > 0$ we can find a number $N$, which may depend on $G$, $H$, and $r$, such that any collection of $N$ distinct left cosets of $H$ in $G$ $\{g_1 H, ..., g_N H\}$ contains a pair $g_i H, g_j H$ which is separated by distance at least $r$, or $d(g_i H, g_j H) \geq r$.

The original definition in [14] does not mention the word metric on $G$ specifically but rather refers to an arbitrary proper left-invariant metric. It turns out that the choice of proper left-invariant metric on $G$ does not matter, since if $d_1$ and $d_2$ are two such metrics on $G$, the subgroup $H \subseteq G$ has bounded packing with respect to $d_1$ if and only if it has bounded packing with respect to $d_2$. See Lemma 2.2 in [14].

It follows immediately from Definition 3.2.1 that any finite index subgroup of $G$ has bounded packing in $G$ since there are no collections of $N$ distinct cosets for $N > |G : H|$. Let us record this observation as:

**Lemma 3.2.1.** *Any finite index subgroup $K$ of a finitely generated countable group $G$ has bounded packing in $G$.*

Our next Lemma concerns the behavior of the property of bounded packing under passing to subgroups or supergroups:

**Lemma 3.2.2.** *Suppose that $H \leq K \leq G$, and $G$ is a finitely generated group.*

*1. If $H$ has bounded packing in $G$, then $H$ has bounded packing in $K$.*

*2. If $H$ has bounded packing in $K$ and $K$ has bounded packing in $G$, then $H$ has bounded packing in $G$.*

*Proof.* See Lemma 2.4 in [14]. □

An important property of bounded packing is that it is preserved by replacing the original subgroup $H$ with a different subgroup $K \subseteq G$ which either contains it or is contained in it with finite index:

**Proposition 3.2.3.** *Let $G$ be a finitely generated group.*

*1. If $H \leq K \leq G$ and $|K : H| < \infty$, then $H$ has bounded packing in $G$ if and only if $K$ has bounded packing in $G$.*

*2. If $H \leq K \leq G$ and $|G : K| < \infty$, then $H$ has bounded packing in $K$ if and only if $H$ has bounded packing in $G$.*

*3. If $H, K \leq G$ and $|G : K| < \infty$, then $H \cap K$ has bounded packing in $K$ if and only if $H$ has bounded packing in $G$.*

*Proof.* See Proposition 2.5 in [14]. □

Part (1) of the proposition has the following immediate corollary:

**Corollary 3.2.4.** *Any finite subgroup of the finitely generated group $G$ has bounded packing.*

*Proof.* It suffices to show that the trivial group $\{1\}$ has bounded packing in $G$, but that is an easy consequence of the properness of the word metric on $G$. □

The following two results summarize the elementary properties of bounded packing. For us, Lemma 3.2.5 will be especially important.

**Lemma 3.2.5.** *Let $1 \longrightarrow N \longrightarrow G \longrightarrow Q \longrightarrow 1$ be a short exact sequence of groups where $G$ is finitely generated. Let $H$ be a subgroup of $G$ which projects to the subgroup* $\overline{H}$ of $Q$. Then $\overline{H}$ has bounded packing in $Q$ if and only if $HN$ has bounded packing in $G$.

*Proof.* See Lemma 2.8 in [14]. □

**Corollary 3.2.6.** *Every normal subgroup $N$ of a finitely generated group $G$ has bounded packing.*

*Proof.* Follows immediately from Lemma 3.2.5 by taking $H = \{1\}$. □

**Theorem 3.2.7.** *If $G$ is a finitely generated virtually nilpotent group, then each subgroup of $G$ has bounded packing in $G$.*

*Proof.* The proof is by induction on the length of the lower central series of $G$ and relies mostly on Lemma 3.2.5. See Theorem 2.12 in [14]. □

The main result of [14] is about bounded packing in *relatively hyperbolic* groups. First, we recall the definition of a relatively hyperbolic group following [14]:

**Definition 3.2.2.** (Relatively hyperbolic group) Let $G$ be a finite generated group and let $\mathbb{P}$ be a finite collection of subgroups of $G$. Suppose that $G$ acts on a $\delta$-hyperbolic graph $\Gamma$ with finite edge stabilizers and finitely many orbits of edges, and suppose further that for each $n \in \mathbb{N}$, each edge in $\Gamma$ is contained in only finitely many circuits of length $n$. If $\mathbb{P}$ is a set of representatives of all the conjugacy classes of infinite vertex stabilizers, then the pair $(G, \mathbb{P})$ is called a *relatively hyperbolic* group. The subgroups in the collection $\mathbb{P}$ are called *peripheral subgroups* and their left cosets are called *peripheral cosets*.

For an equivalent definition and an explanation of the terminology used in the definition above, refer to the extremely readable account by Farb [8].

**Theorem 3.2.8.** *Let $(G, \mathbb{P})$ be a relatively hyperbolic group, with a finite generating set $\mathbb{S}$, and let $H$ be a $\nu$-relatively quasiconvex subgroup of $G$. Suppose that for each* peripheral subgroup $P \in \mathbb{P}$ and each $g \in G$ the intersection $P \cap gHg^{-1}$ has bounded packing in $P$. Then, $H$ has bounded packing in $G$.

*Proof.* See Theorem 8.10 in [14]. □

We do not provide the definition of *relative quasiconvexity* since we shall not make use of Theorem 3.2.8, we refer the reader to the original article [14] instead.

In light of Theorem 3.2.7 one can ask the question whether every subgroup of a solvable group has bounded packing. While an affirmative answer to this question is unlikely, one has the following conjecture promoted in [14], which has recently been settled, see [33] and [28]:

**Conjecture** Every subgroup of a virtually polycyclic group $P$ has the bounded packing property.

The conjecture, in its entirety, was settled in [33] which was announced simultaneously with our proof of the special case which we describe in Section 3.6 below. While Yang's approach does answer the question in the affirmative, it relies heavily on a classical algebraic result by Malcev which states that polycyclic groups are *subgroup separable*. In doing so, however, it gives no insight into the beautiful, and as of yet poorly understood geometry of these groups.

A group $G$ is called *subgroup separable* if each finitely generated subgroup $H \subseteq G$ is an intersection of finite index subgroups of $G$. In [33], Yang approaches the question of bounded packing by establishing that every separable subgroup has bounded packing in its ambient group then quoting Malcev's theorem. By contrast, our approach to solving this conjecture was much more geometric.

## 3.3   Polycyclic Groups and the Groups $P_{n,\phi}$

First, let us define the main object of study:

**Definition 3.3.1.** A group $G$ is called *polycyclic* if there exist subgroups $G_0 = \{1\} \subseteq G_1 \subseteq \cdots \subseteq G_{n-1} \subseteq G_n = G$ with $G_{i-1}$ normal in $G_i$, such that $G_i/G_{i-1}$ is a cyclic group.

An easier to understand from geometrical perspective example of a polycyclic group is the semidirect product $\mathbb{Z}^n \rtimes_\phi \mathbb{Z}$, $\phi \in Aut(\mathbb{Z}^n)$. As a set $\mathbb{Z}^n \rtimes_\phi \mathbb{Z}$ is just $\mathbb{Z}^{n+1}$, however $\mathbb{Z}^n \rtimes_\phi \mathbb{Z}$ is in general far from being quasi-isometric to $\mathbb{E}^{n+1}$.

Recall that the semidirect product $\mathbb{Z}^n \rtimes_\phi \mathbb{Z}$ has the following presentation:

$$\mathbb{Z}^n \rtimes_\phi \mathbb{Z} = \langle \mathbb{Z}^n, t : tat^{-1} = \phi(a), a \in \mathbb{Z}^n \rangle,$$

where $t$ denotes a generator for the $\mathbb{Z}$ factor on the right.

As a matter of convenience, we shall denote the semidirect product of $\mathbb{Z}^n$ with $\mathbb{Z}$ by $P_{n,\phi}$.

The free abelian group $\mathbb{Z}^n$ is a normal subgroup of $P_{n,\phi}$ and the quotient $P_{n,\phi}/\mathbb{Z}^n$ is isomorphic to $\mathbb{Z}$. Hence, we can write $P_{n,\phi}$ as a disjoint union of the left cosets of $\mathbb{Z}^n$ indexed by $\mathbb{Z}$: $P_{n,\phi} = \coprod_{i \in \mathbb{Z}} t^i \mathbb{Z}^n$. We therefore conclude that we can construct the Cayley graph of $P_{n,\phi}$ in the following way: we take a union of $\mathbb{Z}$ copies of the Cayley graph of $\mathbb{Z}^n$ and insert an extra edge at each vertex joining $t^i a \in t^i \mathbb{Z}^n$ to $t^{i+1}\phi^{-1}(a) \in t^{i+1}\mathbb{Z}^n$, which corresponds to multiplication on the right by the generator $t$.

In Section 3.6 below, we show that in a certain sense $P_{n,\phi}$ has only one interesting subgroup from the point of view of bounded packing, namely the cyclic subgroup $H = \langle t \rangle$. Therefore, our first goal will be to show that $H$ has bounded packing in $P_{n,\phi}$.

It is clear from the description of the Cayley graph of $P_{n,\phi}$ that each left coset of $H$ is a transversal to the left cosets of $\mathbb{Z}^n$. On the other hand, the obvious injection $\iota : \mathbb{Z}^n \to t^i\mathbb{Z}^n \subseteq P_{n,\phi}$ induces a graph monomorphism $\iota : \Gamma(\mathbb{Z}^n) \hookrightarrow \Gamma(P_{n,\phi})$. Therefore, in addition to the metric on $t^i\mathbb{Z}^n \subseteq P_{n,\phi}$ coming from the restriction of the word metric $d$ of $P_{n,\phi}$, we have a metric $d_{t^i\mathbb{Z}^n}$ which is the push-forward of the word metric on $\mathbb{Z}^n$. Our first result shows that given any $D > 0$, we can find a $N(D) \in \mathbb{N}$, such that any collection of $N(D)$ distinct left cosets of $H$ will contain a pair $aH, bH$ such that $d_{t^i\mathbb{Z}^n}(at^i, bt^i) > D$ for all $i \in \mathbb{Z}$, which in turn implies $d(aH, bH) > \dfrac{D}{2}$. Then, having shown bounded packing for $H$, we prove bounded packing for all of the subgroups of $P_{n,\phi}$.

In order to show bounded packing of $H$ in $P_{n,\phi}$ as outlined above, we need to analyze the action of $\mathbb{Z}$ on $\mathbb{Z}^n$ via $\phi \in Aut(\mathbb{Z}^n)$. This requires us to present some basic results on the dynamics of linear maps.

## 3.4   The Dynamics of Linear Maps

The main reference for this section is the classical reference on dynamical systems by Hasselblatt and Katok [12]. Their book, however, offers much more than we need, in fact, most of the material we require is contained in Section 1.2 of [12].

Throughout this section, let $\phi : \mathbb{R}^n \to \mathbb{R}^n$ be a linear map. It is well known that the algebraic behavior of a linear map $\phi$ is, to a large extent, determined by its *eigenvalues*. Just for the sake of completeness, we recall that $\lambda \in \mathbb{C}$ is an *eigenvalue* for $\phi$, considered as a linear map $\mathbb{C}^n \to \mathbb{C}^n$, if there exists $v \in \mathbb{C}^n$ such that $\phi(v) = \lambda v$. The set of eigenvalues of $\phi$ is called the spectrum of $\phi$, and the absolute value of the largest eigenvalue of $\phi$, denoted by $\rho(\phi)$, is called the *spectral radius* of $\phi$. In the case of bounded linear maps between Banach spaces, see the discussion below, the spectrum need not be finite. It is, however, always a closed bounded and nonempty subset of $\mathbb{C}$, see for example [16]. For our purposes, we shall only need to consider finite spectra, in particular the spectra of linear maps on $\mathbb{R}^n$. We denote the spectrum of a linear map $\phi$ by $Spec(\phi)$.

If we fix a basis for $\mathbb{R}^n$, such as the standard basis, we can represent the linear map $\phi$ by a matrix $A \in M_n(\mathbb{R})$. To avoid repetition, let us fix a basis once and for all: In what follows, $\{e_i\}_{1 \leq i \leq n}$ will denote the standard basis for $\mathbb{Z}^n$, $\mathbb{R}^n$, and $\mathbb{C}^n$, and the matrix representation of $\phi$ will always be assumed to be with respect to the standard basis unless otherwise stated. We shall use $\phi$ and its matrix representation $A$ interchangeably, whenever there is little chance of confusion.

We would like to understand the behavior of the *iterates* of $\phi$. The linear case is of fundamental importance in the theory of dynamical systems because the derivative functor $D$ establishes a connection between the diffeomorphism $f$ of a Riemannian manifold $M$ and linear maps between finite dimensional vector spaces.

**Definition 3.4.1.** The $n$-fold composition of $\phi$ with itself, $\phi \circ \phi \circ \ldots \circ \phi$ is called the $n$-th *iterate* of $\phi$ and is denoted by $\phi^n$.

More precisely, we are interested in the sets $\mathcal{O}_v^+ = \{\phi^i(v) : i \in \mathbb{N}\}$, for $v \in \mathbb{R}^n$, and $\mathcal{O}_v = \{\phi^i(v) : i \in \mathbb{Z}\}$, where we define $\phi^{-i} = (\phi^{-1})^n$ for $i > 0$, whenever $\phi$ is invertible. The former set is called the *forward orbit* of $v$ under $\phi$, and the latter, simply the *orbit* of $v$ under $\phi$. In the case when $\phi$ is invertible, the forward orbit of $\phi^{-1}$ is called the *backward orbit* of $\phi$. It turns out that the eigenvalues of $\phi$ also determine the asymptotic behavior of the orbits of $\phi$, as the following proposition demonstrates:

**Proposition 3.4.1.** *For every $\delta > 0$, there exists a norm on $\mathbb{R}^n$ such that* $\sup_{||v||=1} ||\phi(v)|| < \rho(\phi) + \delta$.

*Proof.* See Proposition 1.2.2 in [12]. $\qquad\qquad\square$

It seems like Proposition 3.4.1 is too vague to be useful, since it only establishes the bound on the norm of $\phi$ for **some** norm on $\mathbb{R}^n$. However, it is a well know fact that all norms on a finite dimensional normed space are equivalent in the following strong sense:

**Proposition 3.4.2.** *Let $\|\cdot\|_1$ and $\|\cdot\|_2$ be two norms on the finite dimensional vector space $V$, then there exists $C > 0$ such that $\frac{1}{C} \|v\|_1 \leq \|v\|_2 \leq C \|v\|_1$, for all $v \in V$.*

*Proof.* Exercise. □

Before we state the interesting for us corollary to Proposition 3.4.1, let us establish some notation.

The space of all linear maps on a normed vector space $V$ will be denoted by $\mathcal{L}(V)$. It is easily seen that $\mathcal{L}(V)$ is itself a vector space under the usual rules for working with linear maps. Given $\phi \in \mathcal{L}(V)$, the positive real number $\sup_{\|v\|=1} \|\phi(v)\|$ is called the *norm* of the linear map $\phi$. Unfortunately, however, it does not define a norm on all of $\mathcal{L}(V)$ since we can have a linear map on an infinite dimensional vector space $V$ for which this number is not finite. We get around this by defining $\mathcal{B}(V) = \{\phi : V \to V : \|\phi\| < \infty\}$. This subspace of $\mathcal{L}(V)$ is called the *space of bounded linear maps on $V$*, and with the supremum norm defined above, $\mathcal{B}(V)$ becomes a Banach space. A basic result of functional analysis asserts that a linear map is bounded if and only if it is continuous. Actually, $\mathcal{B}(V)$ has more structure than just that of a Banach space, namely it possesses the structure of a *Banach algebra* where multiplication is given by composition of linear maps. It is an easy exercise to show that $|\phi| = \sup_{v \in V} \frac{\|\phi(v)\|}{\|v\|}$, in other words, $|\phi| \, \|v\|$ is an upper bound for $\|\phi(v)\|$ for **all** $v \in V$. Even though we are interested only in the finite dimensional case, we have presented these basics on Banach spaces in order to put the present discussion into proper context. We are now prepared to state the promised corollary to Proposition 3.4.1, which follows immediately in view of Proposition 3.4.2:

**Corollary 3.4.3.** *If $\phi \in \mathcal{L}(\mathbb{R}^n)$ is a bounded linear map with spectral radius $\rho(\phi) < 1$, then $\lim_{i \to \infty} \phi^i(v) = \boldsymbol{0}$, for all $v \in \mathbb{R}^n$. Equivalently, the only accumulation point of $\mathcal{O}^+(v)$ in $\mathbb{R}^n$ is $\boldsymbol{0}$. If in addition $\phi$ is also invertible, then the backward orbit of $\phi$ diverges to infinity.*

*Proof.* See Corollary 1.2.4 in [12]. □

Even though we are not interested in the speed of convergence, the proof of the corollary shows that it is exponential. Let us now consider the following easy to visualize example:

**Example 3.4.4.** Let $\phi : \mathbb{R}^2 \to \mathbb{R}^2$ be the map given by the $2 \times 2$ matrix $\begin{pmatrix} \lambda & 0 \\ 0 & \lambda \end{pmatrix}$, where $0 < \lambda < 1$. The map $\phi$ act by multiplication by $\lambda < 1$ and for every $v \in \mathbb{R}^2$, $\phi^i(v)$ converges to $\mathbf{0}$ exponentially as $i \to \infty$.

We are now ready to present our motivating example for much of the discussion that follows:

**Example 3.4.5.** Let $0 < \mu < 1 < \lambda$ and consider the map $\phi : \mathbb{R}^2 \to \mathbb{R}^2$ given by $A = \begin{pmatrix} \mu & 0 \\ 0 & \lambda \end{pmatrix}$. The matrix $A$ has an eigenbasis $\{e_1, e_2\}$ with eigenvalues $\mu$ and $\lambda$ respectively. Let us define $E_- = \mathbb{R}e_1$ and $E_+ = \mathbb{R}e_2$. The forward orbit of any point in $E_-$ approaches the origin at an exponential speed, and the same is true for the backward orbit of any point in $E_+$.

Now, we proceed to analyze the behavior of the orbit of a point $(x_0, y_0) \in \mathbb{R}^2 - (E_- \cup E_+)$ and to that end, without loss of generality, we assume that $x_0, y_0 > 0$. We compute $\phi^i(x_0, y_0) = (\mu^i x_0, \lambda^i y_0)$, hence $\phi^i(x_0, y_0)$ lies on the curve $\mathcal{F}_{(x_0, y_0)} = \{(x, y) \in \mathbb{R}^2 : xy^c = x_0 y_0^c\}$, where the constant $c > 0$ is defined by $\mu = \lambda^{-c}$. Further, it is easily shown that $\phi$ leaves the curves $\mathcal{F}_{(x_0, y_0)}$ invariant. In the special case where $\mu = \lambda^{-1}$, the curves $\mathcal{F}_{(x_0, y_0)}$ give a foliation of the plane by hyperbolae.

Example 3.4.5 motivates the following definition:

**Definition 3.4.2.** (Hyperbolic map) A linear map $\phi : \mathbb{R}^n \to \mathbb{R}^n$ is called *hyperbolic*, if the intersection of $Spec(\phi)$ with the unit circle in $\mathbb{C}$ is empty, or in other words, $\phi$ has no absolute value 1 eigenvalues.

Next, we recast some familiar notions from linear algebra from the point of view of dynamics. Recall that to each $\lambda \in Spec(\phi) \cap \mathbb{R}$, for a given $\phi \in \mathcal{L}(\mathbb{R}^n)$, one associates its *root space* $\bigcup_k Ker(\phi - \lambda I)^k$ denoted by $E_\lambda$. Similarly, to each pair of complex conjugates $\lambda, \overline{\lambda} \in Spec(\phi)$, we define $E_{\lambda, \overline{\lambda}}$ to be the intersection of the root space of the complexified map $\phi \in \mathcal{L}(\mathbb{C}^n)$ with $\mathbb{R}^n$. By analogy with Example 3.4.5, we set

$$E_- = \bigoplus_{|\lambda|<1} E_\lambda \oplus \bigoplus_{|\lambda|<1} E_{\lambda,\overline{\lambda}}, \tag{3.1}$$

$$E_+ = \bigoplus_{|\lambda|>1} E_\lambda \oplus \bigoplus_{|\lambda|>1} E_{\lambda,\overline{\lambda}}, \tag{3.2}$$

and

$$E_0 = E_{-1} \oplus E_1 \oplus \bigoplus_{|\lambda|=1} E_{\lambda,\overline{\lambda}}. \tag{3.3}$$

Clearly then, $\mathbb{R}^n = E_- \oplus E_+ \oplus E_0$, and also $\phi$ is a hyperbolic map if and only if $E_0 = \{\mathbf{0}\}$.

Everything we need from [12] is contained in the following:

**Proposition 3.4.6.** *Let $\phi \in \mathcal{L}(\mathbb{R}^n)$ be a hyperbolic map. Then:*

1. *For every $v \in E_-$, the positive iterates $\phi^i(v)$ converge to the origin with exponential speed as $i \to \infty$, and if $\phi$ is invertible, then the negative iterates $\phi^i(v)$ go to infinity with exponential speed as $i \to -\infty$.*

2. *For every $v \in E_+$, the positive iterates of $v$ go to infinity exponentially and if $\phi$ is invertible, then the negative iterates converge exponentially to the origin.*

3. *For every $v \in \mathbb{R}^n - (E_- \cup E_+)$, the iterates $\phi^i(v)$ converge to infinity exponentially as $i \to \pm\infty$.*

*Proof.* See Proposition 1.2.8 in [12].   $\square$

## 3.5   The Diagonal Action of $\mathbb{Z}$ on $\mathbb{Z}^n$

This section contains a couple of technical lemmas which analyze the geometry of the space of left cosets of $H$ and lays a foundation for the work in Sections 3.6 - 3.7.

Any automorphism of $\mathbb{Z}^n$ extends to an automorphism of $\mathbb{R}^n$. On the other hand, an automorphism of $\mathbb{R}^n$ which preserves the integer lattice $\mathbb{Z}^n \subseteq \mathbb{R}^n$ restricts to an automorphism of $\mathbb{Z}^n$ if and only if it has determinant equal to $\pm 1$. Thus, we identify $Aut(\mathbb{Z}^n)$ with the group of $n \times n$ integer matrices of determinant $\pm 1$. This group is denoted by $GL_n(\mathbb{Z})$ and as we noted $GL_n(\mathbb{Z}) \subseteq GL_n(\mathbb{R})$.

Now, given $\phi \in Aut(\mathbb{Z}^n)$ we turn our attention to the action of $\mathbb{Z}$ on $\mathbb{Z}^n$ induced by the homomorphism $\Phi : \mathbb{Z} \to Aut(\mathbb{Z}^n)$ given by $\Phi : 1 \mapsto \phi$. If the automorphism $\phi$ is diagonalizable over $\mathbb{R}$, we shall say that $\mathbb{Z}$ *acts diagonally* on $\mathbb{Z}^n$. Our plan of attack on the problem of bounded packing of $H$ in $P_{n,\phi}$, as outlined at the closing of Section 3.3, is to show that if one takes a large enough collection of distinct left cosets $a_1 H, ..., a_N H$, $a_i \in \mathbb{Z}^n$, then one is guaranteed that some pair will stay far apart in each left coset of $\mathbb{Z}^n$ in the metric $d_{t^i \mathbb{Z}^n}$. It turns out that this is intimately connected with the geometry of the orbits of the $\mathbb{Z}$ action on $\mathbb{Z}^n$ via $\phi$.

We begin by showing that given a diagonal action of $\mathbb{Z}$ on $\mathbb{Z}^n$ by a hyperbolic automorphism, two translates of orbits of this action intersect in at most two points. The reason for considering diagonal actions is simple, if all of the eigenvalues of $\phi$ are real and positive, the action of $\mathbb{Z}$ on $\mathbb{Z}^n$ extends to a *flow* on $\mathbb{R}^n$ and we can treat the orbits of the $\mathbb{Z}$-action as differentiable curves in $\mathbb{R}^n$.

**Lemma 3.5.1.** *Let $\mathbb{Z}$ act diagonally on $\mathbb{Z}^n$ via the hyperbolic map $\phi \in Aut(\mathbb{Z}^n)$ all of whose eigenvalues are positive. If $\mathcal{O}_z$ and $\mathcal{O}_w$ are the orbits of the points $z, w \in \mathbb{Z}^n - \{0\}$, then $|\mathcal{O}_z \cap (a + \mathcal{O}_w)| \leq 2$ for any non-zero $a \in \mathbb{Z}^n$.*

*Proof.* Let $\{v_1, ..., v_n\}$ be an eigenbasis for $\phi$. For $z \in \mathbb{Z}^n$, write $z = z_1 v_1 + ... z_n v_n$, $a = a_1 v_1 + ... + a_n v_n$. Then, $\mathcal{O}_z$ is contained in the image of the curve $t \to (\lambda_1^t z_1, ..., \lambda_n^t z_n)$.

Now, given $z, w \in \mathbb{Z}^n$, set $\alpha(t) = (\lambda_1^t z_1, ..., \lambda_n^t z_n)$ and $\beta(t) = (\lambda_1^t w_1 + a_1, ..., \lambda_n^t w_n + a_n)$, and note that $\mathcal{O}_z \cap (a + \mathcal{O}_w) \subseteq Im(\alpha) \cap Im(\beta)$. We show that there are at most two solutions $(t, t')$ to the equation $\alpha(t) = \beta(t')$. Since $a \neq 0$, at least two coordinates of $a$ with respect to the chosen eigenbasis for $\phi$ are non-zero (otherwise $a$ would lie in the intersection of a 1-dimensional eigenspace of $\phi$ with $\mathbb{Z}^n$, which is impossible). By reordering the basis, if necessary, we may assume that $a_1 \neq 0$ and $a_2 \neq 0$. We consider the following cases:

1. $z_1 = 0$ and $w_1 \neq 0$: Considering the equation $\lambda_1^t z_1 = \lambda_1^{t'} w_1 + a_1$, we see that there is at most one solution for $t'$. Now, we consider the subcases:

   - $w_2 = 0$: Consider the equation $\lambda_2^t z_2 = a_2$. If $z_2 = 0$, this equation has no solution. If $z_2 \neq 0$, there exists at most one solution for $t$.

   - $w_2 \neq 0$ and $z_2 \neq 0$: Considering $\lambda_2^t z_2 = \lambda_2^{t'} w_2 + a_2$, with the value for $t'$ we found above, we again conclude that there is at most one solution for $t$.

   - $w_2 \neq 0$ and $z_2 = 0$: In this case, since $z_1 = z_2 = 0$, we must have some $z_i \neq 0$ since $z \neq 0$. Assume that $z_3 \neq 0$, and consider $\lambda_3^t z_3 = \lambda_3^{t'} w_3 + a_3$. With the value for $t'$, there is again at most one solution for $t$.

2. $z_1 \neq 0$, $w_1 = 0$, $z_2 \neq 0$, $w_2 \neq 0$: Consider $\lambda_1^t z_1 = a_1$. We conclude that there is at most one solution for $t$. With this value for $t$ we solve $\lambda_2^t z_2 = \lambda_2^{t'} w_2 + a_2$, and we find at most one solution for $t'$.

3. $z_1 \neq 0$, $w_1 \neq 0$, $z_2 \neq 0$, $w_2 \neq 0$: Consider the projections of the curves $\alpha$ and $\beta$ on the plane spanned by $\{v_1, v_2\}$. These projections are the curves $t \to (\lambda_1^t z_1, \lambda_2^t z_2)$ and
$t \to (\lambda_1^t w_1 + a_1, \lambda_2^t w_2 + a_2)$ respectively. Setting $t' = t + k$, $\lambda_2 = \lambda_1^p$, $\lambda_1^t = x$, $\lambda_1^k = y$, we obtain the system:

$$xz_1 = xyw_1 + a_1$$
$$x^p z_2 = x^p y^p w_2 + a_2$$

Solving for $xy$ and substituting, we get:

$$x^p z_2 = \left(\frac{xz_1 - a_1}{w_1}\right)^p w_2 + a_2.$$

We now show that this equation has at most 2 solutions for $x$. Consider the function $f(x) = x^p z_2 - \left(\frac{xz_1 - a_1}{w_1}\right)^p w_2 - a_2$. Then, $f'(x) = p\left[z_2 x^{p-1} - \frac{z_1 w_2}{w_1}\left(\frac{xz_1 - a_1}{w_1}\right)^{p-1}\right]$. Note that $p \neq 0$, as no eigenvalue of $\phi$ is equal to 1. Note also that we may always arrange to have $p \neq 1$, since if $p = 1$, then $\lambda_1 = \lambda_2$, in which case we look for $a_i$, with $i \neq 1, 2$, such that $a_i \neq 0$, and $\lambda_i \neq \lambda_1, \lambda_2$, and replace $v_2$ with $v_i$. If no such $a_i$ exists, we conclude that $a$ lies in the intersection of $\mathbb{Z}^n$ and an eigenspace of $\phi$, which is impossible. As $a_1 \neq 0$, the only possible zero of $f'$ must be located at $x = \frac{a_1}{z_1 - w_1\left(\frac{z_2 w_1}{z_1 w_2}\right)^{\frac{1}{p-1}}}$. Since $f$ has at most one more zero than $f'$, we conclude that $f$ has at most two zeros. Therefore, we have at most two solutions for $t$, and solving for $y$ concludes this final case.

It is an easy exercise for the reader to convince himself that this exhausts all possible cases up to permutation of the variables involved. $\square$

Even though Lemma 3.5 was stated to handle the case where only one of the orbits was translated by $a \in \mathbb{Z}^n$, it takes no work at all to show that for any $z, w, a, b \in \mathbb{Z}^n$, we still have $|(a + \mathcal{O}_z) \cap (b + \mathcal{O}_w)| \leq 2$. Now, using this counting argument, we prove the following interesting property of the diagonal action of $\mathbb{Z}$. Consider the induced action of $\phi$ on the set $\mathcal{F}(\mathbb{Z}^n)$ of finite subsets of $\mathbb{Z}^n$. We show that given $D > 0$, there exist $N(D) \in \mathbb{N}$, such that any $\mathcal{C} \in \mathcal{F}(\mathbb{Z}^n)$ with $|\mathcal{C}| > N(D)$ is guaranteed to contain a pair of points whose images in every $\phi^k(\mathcal{C})$ are separated by a distance of at least $D$. In other words, there exist $a, b \in \mathcal{C}$ such that $d(\phi^k(a), \phi^k(b)) > D$ for all $k \in \mathbb{Z}$.

**Lemma 3.5.2.** *Let $\mathbb{Z}$ act diagonally on $\mathbb{Z}^n$ via the hyperbolic automorphism $\phi$ all of whose eigenvalues are positive. Then given $D > 0$, there exists $N > 0$ such that for any collection of $m > N$ distinct points $a_1, a_2, ..., a_m \in \mathbb{Z}^n$, there exist $i, j$ such that $d(\phi^k(a_i), \phi^k(a_j)) > D$ for all $k \in \mathbb{Z}$.*

*Proof.* Let $S$ be the set $\overline{B}(0, D) \cap \mathbb{Z}^n$. If $\left\| \phi^k(z_1) - \phi^k(z_2) \right\| \leq D$ for some $z_1, z_2 \in \mathbb{Z}^n$ and $k \in \mathbb{Z}$, then $z_1 - z_2$ is in the orbit $\mathcal{O}_z$ of some $z \in S$ under $\phi$. Hence, $z_2 \in z_1 + \bigcup_{z \in S} \mathcal{O}_z$. Now, if we have $\left\| \phi^k(a_i) - \phi^k(a_j) \right\| \leq D$ for all $1 \leq i, j \leq m$, then $a_m \in \bigcap_{1 \leq l \leq m-1}(a_l + \bigcup_{z \in S} \mathcal{O}_z)$. To finish the proof, we observe that Lemma 3.5.1 shows that $\left| \bigcap_{1 \leq l \leq m-1}(a_l + \bigcup_{z \in S} \mathcal{O}_z) \right| \leq 2 \left| S \right|^2 \sim 2 \left( Vol(\overline{B}(0, D)) \right)^2$. $\qquad \square$

## 3.6   Bounded Packing in $P_{n,\phi}$

We now use the dynamical properties of the diagonal $\mathbb{Z}$-action on $\mathcal{F}(\mathbb{Z}^n)$ to prove bounded packing in $\mathbb{Z}^n \rtimes \mathbb{Z}$.

**Lemma 3.6.1.** *Let $\mathbb{Z}$ act diagonally on $\mathbb{Z}^n$ via the hyperbolic automorphism $\phi$, all of whose eigenvalues are positive. Let $P_{n,\phi} = \mathbb{Z}^n \rtimes_\phi \mathbb{Z}$ and let $t$ denotes a generator for the copy of $\mathbb{Z}$ in $G$ on the right, then the subgroup $H = \langle t \rangle$ has bounded packing in $P_{n,\phi}$.*

*Proof.* Let $d$ denote the word metric on $G$, and let $d_{t^i \mathbb{Z}^n}$ denote the inherited Euclidean metric on $t^i \mathbb{Z}^n$ described in Section 3.3. We show that given $R > 0$, there exists $D > 0$ such that if $x, y \in t^i \mathbb{Z}^n$ and $d_{t^i \mathbb{Z}^n}(x, y) > D$, then $d(x, y) > R$. After translating, we need only show that if $d_{\mathbb{Z}^n}(e, y) > D$, then $d(e, y) > R$. Supposing to the contrary, that given a fixed $R > 0$, for any $D > 0$ there exists $y \in \mathbb{Z}^n$ with $d_{\mathbb{Z}^n}(e, y) > D$ and $d(e, y) < R$, we construct a sequence $\{y_k\}$ with $d_{\mathbb{Z}^n}(e, y_k) > k$ and $d(e, y_k) < R$. Because the metric $d$ is proper, this is a contradiction. Now, since $\mathbb{R}^n$ is quasi-isometric to the lattice $\mathbb{Z}^n$, Lemma 3.5.2 shows that given $D > 0$, there exists $N \in \mathbb{N}$ such that any collection of $m$ distinct cosets of $H$, with $m > N$, contains a pair, say $aH$, $bH$, such that in any coset $t^i\mathbb{Z}^n$, we have $d_{t^i\mathbb{Z}^n}(at^i, bt^i) > D$ for all $i \in \mathbb{Z}$. This along with the discussion above shows that given $D > 0$ we can find $N' \in \mathbb{N}$ such that any collection of $m$ distinct cosets, with $m > N'$, contains a pair, $aH$, $bH$, such that $d(at^i, bt^i) > D$ for all $i \in \mathbb{Z}$. Finally, we show that this implies that $d(aH, bH) \geq \dfrac{D}{2}$. To this end, let $at^i \in aH$, and $bt^j \in bH$. We have two cases: either $|i - j| \geq \dfrac{D}{2}$, or $|i - j| < \dfrac{D}{2}$. In the first case, it is clear that $d(at^i, bt^j) \geq |i - j| \geq \dfrac{D}{2}$. In the second, we have: $d(at^i, bt^j) \geq |d(at^i, at^j) - d(at^j, bt^j)| > \left| \dfrac{D}{2} - D \right| = \dfrac{D}{2}$. $\qquad\square$

Having established bounded packing of $H$ in $P_{n,\phi}$ where the diagonal $\mathbb{Z}$-action on $\mathbb{Z}^n$ is via a hyperbolic automorphism, we now work to relax the hyperbolicity assumption.

**Lemma 3.6.2.** *Let $\mathbb{Z}$ act diagonally on $\mathbb{Z}^n$ via the automorphism $\phi$ all of whose eigenvalues are positive. If $P_{n,\phi} = \mathbb{Z}^n \rtimes_\phi \mathbb{Z}$ and if $t$ denotes a generator for the right copy of $\mathbb{Z}$ in $P_{n,\phi}$, then the subgroup $H = \langle t \rangle$ has bounded packing in $P_{n,\phi}$.*

*Proof.* The automorphism $\phi$ extends to an automorphism $\phi : \mathbb{R}^n \mapsto \mathbb{R}^n$ which has determinant equal to $\pm 1$. We have $\mathbb{R}^n = \mathcal{H} \oplus E_1$, where $\mathcal{H} = E_- \oplus E_+$ is the sum of the expanding and contracting subspaces of $\phi$, and $E_1$ is the eigenspace for the eigenvalue $\lambda = 1$ of $\phi$. Note that since both $\mathcal{H}$ and $E_1$ are $\phi$-stable, so is $Z_\mathcal{H} = \mathbb{Z}^n \cap \mathcal{H}$, so $Z_\mathcal{H} \lhd P_{n,\phi}$. Because $Z_\mathcal{H} = \mathbb{Z}^n \cap \mathcal{H}$ is a free abelian group on which $\phi$ acts as a hyperbolic automorphism, $\langle t \rangle$ has bounded packing in $Z_\mathcal{H} \rtimes_\phi \mathbb{Z} = Z_\mathcal{H} \langle t \rangle \subseteq P_{n,\phi}$ by Lemma 3.6.1. Now, $P_{n,\phi}/Z_\mathcal{H}$ embeds in $\left( \mathbb{R}^n \rtimes_\phi \mathbb{R} \right) / \mathcal{H} \cong (\mathbb{R}^n / \mathcal{H}) \rtimes_{\overline{\phi}} \mathbb{R}$ which is isomorphic to $\mathbb{R}^k$ since $\overline{\phi}$ is the identity. Therefore, $P_{n,\phi}/Z_\mathcal{H}$ is abelian hence it has bounded packing with respect to all of its subgroups. In particular, the image of $\langle t \rangle$ has bounded packing in $P_{n,\phi}/Z_\mathcal{H}$, and therefore $Z_\mathcal{H} \langle t \rangle$ has bounded packing in $P_{n,\phi}$. As $\langle t \rangle \subseteq Z_\mathcal{H} \langle t \rangle \subseteq P_{n,\phi}$, we conclude that $\langle t \rangle$ has bounded packing in $P_{n,\phi}$. $\qquad\square$

An immediate consequence of Lemma 3.6.2 is:

**Corollary 3.6.3.** *With the notation from Lemma 3.6.2, for every $z \in \mathbb{Z}^n$, the subgroup $\langle zx \rangle$ has bounded packing in $P_{n,\phi}$.*

*Proof.* Since for any $z \in \mathbb{Z}^n$, we have:

$zxw\,(zx)^{-1} = z\phi(w)z^{-1} = \phi(w) = xzx^{-1}$, we have an automorphism of $P_{n,\phi}$ which is the identity on $\mathbb{Z}^n$ and which sends $x$ to $zx$. The conclusion now follows immediately from Lemma 3.6.2. □

This corollary shows that any cyclic subgroup of $P_{n,\phi}$ has bounded packing.

We finally have the means to prove the main result of this section:

**Theorem 3.6.4.** *Let $\mathbb{Z}$ act diagonally on $\mathbb{Z}^n$ via $\phi \in Aut(\mathbb{Z}^n)$ all of whose eigenvalues are real, then every subgroup of $P_{n,\phi} = \mathbb{Z}^n \rtimes_\phi \mathbb{Z}$ has bounded packing in $P_{n,\phi}$.*

*Proof.* By passing to a subgroup of finite index, we may assume that all of the eigenvalues of $\phi$ are real and positive. Let $H \subseteq P_{n,\phi}$. By further passing to a finite index subgroup of $P_{n,\phi}$, we may assume that the projection of $H$ to $\mathbb{Z}$ is onto, so that we have the exact sequence $1 \longrightarrow W \longrightarrow H \longrightarrow \mathbb{Z} \longrightarrow 1$, where $W = \mathbb{Z}^n \cap H$. Note that $W$ is normal in $P_{n,\phi}$: let $zt \in H$ be an element of $H$ which projects to a generator of the right-side copy of $\mathbb{Z}$ in $P_{n,\phi}$, where $z \in \mathbb{Z}^n$, then $ztWt^{-1}z^{-1} = W$, so $tWt^{-1} = z^{-1}Wz = W$. Next, by passing to a subgroup of $P_{n,\phi}$ containing $H$ with finite index, we may assume that $\mathbb{Z}^n/W$ is free abelian. Here is how: let $\mathbb{Z}^n/W = F \oplus T$, where $F$ is free abelian, and $T$ is torsion; the preimage of $T$ in $\mathbb{Z}^n$, say $T'$, under the quotient map is $\phi$-stable, $|T' : W| < \infty$, so $T' \rtimes \mathbb{Z}$ is the desired subgroup. Now, consider $G = P_{n,\phi}/W$. Because the diagram

$$
\begin{array}{ccccccccc}
1 & \longrightarrow & \mathbb{Z}^n & \longrightarrow & P_{n,\phi} & \longrightarrow & \mathbb{Z} & \longrightarrow & 1 \\
& & \downarrow & & \downarrow & & & & \\
1 & \longrightarrow & \mathbb{Z}^n/W & \longrightarrow & P_{n,\phi}/W & \longrightarrow & \mathbb{Z} & \longrightarrow & 1
\end{array}
$$

commutes, $G \cong (\mathbb{Z}^n/W) \rtimes_{\overline{\phi}} \mathbb{Z} \cong \mathbb{Z}^m \rtimes_{\overline{\phi}} \mathbb{Z}$, where $\overline{\phi}$ is the induced automorphism on the quotient $\mathbb{Z}^n/W$. Corollary 3.6.3 shows that the image of $\langle zt \rangle$ has bounded packing in $P_{n,\phi}/W$, hence $H = W \langle zt \rangle$ has bounded packing in $P_{n,\phi}$. $\qquad \square$

Note that Theorem 3.6.4 solves the question of bounded packing in polycyclic groups of length $\leq 3$ as the following corollary demonstrates:

**Corollary 3.6.5.** *If $G$ is a polycyclic group of length $\leq 3$, $G$ has bounded packing with respect to all of its subgroups.*

*Proof.* The statement is trivially true if the length of $G$ is equal to 1 or 2, since in both cases $G$ is virtually abelian. If $G$ has length 3, then $G$ contains a subgroup isomorphic to $\mathbb{Z}^2 \rtimes_{\phi} \mathbb{Z}$ of finite index. Unless $\phi$ is hyperbolic, $G$ will be virtually nilpotent. In the case where $\phi$ is hyperbolic, the eigenvalues are real, and the conclusion follows from Theorem 3.6.4. $\qquad \square$

## 3.7 Coset Growth of $H$ in $P_{2,\phi}$

### 3.7.1 Growth of Groups and Coset Growth

**Growth of Groups**

In the course of proving bounded packing in $P_{n,\phi}$, we have developed the tools for studying the growth of a certain metric space. Recall that given a metric measure space $(X, d, \mu)$, one defines the *volume growth* to be the function $f_{x_0}(r) = \mu(B(x_0, r))$, where $B(x_0, r)$ is the open ball of radius $r$ centered at the point $x_0 \in X$. In geometric group theory, the metric space is a finitely generated group $G$ with the word metric, the measure on $G$ is the counting measure, and the basepoint is the identity element $1 \in G$. In other words, the function $f(r)$ counts the number of elements which can be expressed as a word of length at most $r$ in the generators. In order to make the definition independent of the generating set, we introduce the following equivalence relation on the set of non-decreasing functions defined on $\mathbb{R}_{\geq 0}$: $\alpha \sim \beta$ if there exists $C > 0$ such that $\beta\left(\dfrac{r}{C}\right) \leq \alpha(r) \leq \beta\left(Cr\right)$. We now define the *growth rate* of $G$ to be the equivalence class of the growth function $f(r)$.

It is important to note that the growth rate of any finitely generated group is at most exponential. This is due to the fact that the growth rate of any quotient of the finitely generated group $G$ is bounded above by the growth rate of $G$. On the other hand, every finitely generated group is a quotient of a finitely generated free group and these are easily shown to have exponential growth. What else can be said of the growth rate? The free abelian groups $\mathbb{Z}^n$ provide examples of groups of polynomial growth of any order $n$. In [32], Wolf proves that nilpotent groups have polynomial growth. Actually, [32] contains a stronger result which is of interest for us:

**Theorem 3.7.1.** *A polycyclic group is either virtually nilpotent and is thus of polynomial growth, or is not virtually nilpotent and is of exponential growth.*

*Proof.* See [32]. □

Coupling this with a result by Milnor, see [19], Wolf deduces:

**Theorem 3.7.2.** *A finitely generated solvable group is either polycyclic and virtually nilpotent and is thus of polynomial growth, or has no nilpotent subgroup of finite index and is of exponential growth.*

*Proof.* See [32]. □

Much progress has been made in the area of geometric group theory concerned with the growth of groups since the works of Wolf and Milnor. The most notable result is, of course, Gromov's celebrated result which characterizes the finitely generated, virtually nilpotent groups, as precisely the finitely generated groups of polynomial growth:

**Theorem 3.7.3.** *A finitely generated group has polynomial growth if and only if it is virtually nilpotent.*

*Proof.* See [10] or [15]. □

Gromov's original proof of this amazing result relies on the solution of Hilbert's Fifth Problem which, in its original version, aimed to characterize Lie groups as the topological groups which are also topological manifolds, a very powerful result in and of itself. Since Gromov's original proof in [10], simpler proofs, which do not use the solution to Hilbert's Fifth Problem, have appeared. See, for example, the work by Kleiner in [15].

A question which took some time to resolve was whether there existed groups of intermediate growth, or in other words, groups whose growth rate is less than exponential but greater than polynomial. The question was answered in the affirmative by Grigorchuk in [9].

**Coset Growth**

Let $G$ be a finitely generated group and $H$ be any subgroup. Consider the space of left cosets $G/H \subseteq \mathcal{P}(G)$. The word metric $d$ on $G$ induces a distance function on $\mathcal{P}(G)$, see Section 3.1.2. We define the *coset growth* $f_{G,H}$ to be the growth rate of the coset space $G/H$ using the distance function $d$ instead of a true metric. We immediately observe that $f_{G,\{1\}}$ is just the growth rate of $G$, while $f_{G,N}$ is the growth rate of the quotient group $G/N$ whenever $N$ is normal in $G$. We shall also refer to the coset growth of $H$ in $G$ as the growth of the pair $(G, H)$. Having an infinite coset growth of $H$ in $G$ is equivalent to $H$ not having the property of bounded packing in $G$. Our goal in this section is to show that the growth rate of the pair $(P_{2,\phi}, \langle t \rangle)$ is at most exponential, when $P_{2,\phi}$ is not nilpotent. Theorem 3.7.1 shows that the growth rate of $(P_{2,\phi}, \{1\})$ is exponential, whereas the growth rate of $(P_{2,\phi}, \mathbb{Z}^2)$ is polynomial of order 1 as $P_2/\mathbb{Z}^2 \cong \mathbb{Z}$.

In order to bound the growth rate of $(P_{2,\phi}, \langle t \rangle)$ above, we need to examine the distortion of $\mathbb{Z}^2$ in $P_{2,\phi}$. To do this, we make use of the fact that $P_{2,\phi}$ embeds as a uniform lattice of finite covolume in the Lie group *Sol*. In fact, every polycyclic group can be virtually embedded as a lattice in a simply-connected solvable Lie group. However, in the case of $P_{2,\phi}$ the embedding is far less abstract.

### 3.7.2 The *Sol* Geometry

**The World of Geometric Manifolds**

Quite possibly the first mention of the solvable Lie group *Sol* in a mathematician's vocabulary is as one of the eight *Thurston geometries*.

Recall that a *geometry* is a pair $(M, G)$, where $M$ is manifold and $G$ a Lie group such that $G$ acts transitively on $M$ with compact point stabilizers. We say that the geometries $(M_1, G_1)$, $(M_2, G_2)$ are *equivalent* if there exists a diffeomorphism $\phi : M_1 \to M_2$ and an isomorphism $\phi_* : G_1 \to G_2$ such that $\phi(g \cdot x) = \phi_*(g) \cdot \phi(x)$, for all $x \in M_1$ and for all $g \in G_1$. This compact definition captures both the essence of geometry as we understood it from the ancient Greeks and its more contemporary formulation underlying the Erlangen program. Its importance is witnessed by Thurston's famous *geometrization conjecture*, which is now a theorem due to Perelman. Very informally, the geometrization theorem states that a closed orientable 3-manifold can be cut up into pieces which carry a geometry.

There are two well-known algorithms for cutting a 3-manifold into simpler pieces:

**Theorem 3.7.4.** *(Kneser, Milnor) Let $M$ be a closed, orientable 3-manifold. Then, $M$ admits a finite connected sum decomposition $(K_1 \sharp ... \sharp K_p) \sharp (L_1 \sharp ... \sharp L_q) \sharp (\mathbb{S}^2 \times \mathbb{S}^1 \sharp ... \sharp \mathbb{S}^2 \times \mathbb{S}^1)$, where $K_i$ is an irreducible 3-manifold with infinite $\pi_1$ and $L_j$ is an irreducible 3-manifold with finite $\pi_1$, for all $1 \leq i \leq p$, $1 \leq j \leq q$.*

A manifold is said to be *closed* if it is compact and has no boundary. The *boundary*

of a smooth manifold is the subset of $M$ consisting of points which have neighborhoods diffeomorphic to an open neighborhood of the origin in $\mathbb{R}^3_+ = \{(x, y, z) : x, y, z \in \mathbb{R}, z \geq 0\}$. A 3-manifold is called *irreducible* if every embedded 2-sphere is the boundary of an embedded copy of $B^3 = \{\mathbf{x} \in \mathbb{R}^3 : ||\mathbf{x}|| \leq 1\}$. The symbol $\sharp$ in Theorem 3.7.4 denotes the operation of forming a connected sum: $M_1 \sharp M_2 = (M_1 - int B^3) \coprod_{\mathbb{S}^2} (M_2 - int B^3)$. In other words, one excises the interior of a 3-ball from both $M_1$ and $M_2$ then glues the resulting manifolds with boundary along their corresponding copies of $\partial B^3 = \mathbb{S}^2$. A manifold which cannot be written as a connected sum of two manifolds unless one of them is $\mathbb{S}^3$ is called *prime*. In addition to the connected sum decomposition result above, we have the following theorem proved by Jaco-Shalen and independently by Johannson:

**Theorem 3.7.5.** *(Jaco-Shalen, Johannson) Let $M$ be a closed, orientable, irreducible 3-manifold. Then, there is a finite collection of disjoint incompressible tori in $M$ that separate $M$ into a finite collection of compact 3-manifolds with toral boundary each of which is either torus-irreducible or Seifert fibered.*

A Seifert fibered space $X$ is a manifold which can be written as a disjoint union of circles such that each circle has a regular neighborhood diffeomorphic to a *fibered solid torus $D^2 \times I \backslash (\phi(z, 0) \sim \phi(z, 1))$*, where $D^2 = \{z \in \mathbb{C} : ||z|| \leq 1\}$, and $\phi : D^2 \to D^2$ is given by $\phi(z) = e^{2\pi i \frac{p}{q}} z$ for $p, q \in \mathbb{Z}$, see [30].

We are finally ready to state the geometrization theorem:

**Theorem 3.7.6.** *(Perelman) Let $M$ be a closed, oriented 3-manifold. Then, each component obtained by performing the sphere and then the torus decomposition admits a geometric structure.*

For the proof, we refer the reader to the original papers by Perelman: [22], [23], and [24].

Any manifold which admits a geometric structure is covered by one of the eight 3-dimensional simply-connected geometries: $\mathbb{E}^3$, $\mathbb{H}^3$, $\mathbb{S}^3$, $\mathbb{S}^2 \times \mathbb{R}$, $\mathbb{H}^2 \times \mathbb{R}$, $\widetilde{SL_2(\mathbb{R})}$, $Nil$, and $Sol$. An excellent account of the topic can be found in [30].

## Definition and Basic Properties of *Sol*

The Lie group *Sol* is defined to be the semidirect product of $\mathbb{R}^2 \rtimes \mathbb{R}$, where a generator $t \in \mathbb{R}$ acts on $\mathbb{R}^2$ via $\begin{pmatrix} e^t & 0 \\ 0 & e^{-t} \end{pmatrix}$. We identify *Sol* as a set with $\mathbb{R}^3 = \{(x, y, z) : (x, y) \in \mathbb{R}^2, z \in \mathbb{R}\}$ and thereby obtain a chart on the underlying manifold. We shall use this chart to introduce a left-invariant Riemannian metric on *Sol* as follows: First, we metrically identify the normal copy of $\mathbb{R}^2$ with $\mathbb{E}^2$. Next, we transport the metric $ds^2_{z=0} = dx^2 + dy^2 + dz^2$ on $\mathbb{R}^2 \times \{0\}$ to $\mathbb{R}^2 \times \{t\}$ for arbitrary $t \in \mathbb{R}$ using the multiplication map $m_t : Sol \to Sol$ given by $m_t(x, y, z) = (0, 0, t)(x, y, z)$. This means that we define $ds^2_{z=t}$ in such a way that it satisfies $ds^2_{z=t}(Dm_t(\mathbf{v})) = ds^2_{z=0}(\mathbf{v})$ for every $\mathbf{v} \in T_p(\mathbb{R}^2) \subseteq T_p(Sol)$. It is easy to show that the resulting left-invariant Riemannian metric is given by $ds^2 = e^{2z}dx^2 + e^{-2z}dy^2 + dz^2$.

We shall denote the associated path metric on *Sol* by $d_S$. The embedding of $\mathbb{E}^2$ into *Sol* via the map $(x, y) \mapsto (x, y, z)$ endows $\mathbb{R}^2 \subseteq Sol$ with the Euclidean metric which we shall denote by $d_{\mathbb{R}^2}$. In order to produce the promised bound for the growth rate of the pair $(P_2, \langle t \rangle)$, we need the following result which provides a lower bound for the distance in *Sol* between two points in $\mathbb{R}^2$ in terms of their $l^1$ distance in $\mathbb{R}^2$:

**Lemma 3.7.7.** $d_S((x_1, y_1, 0), (x_2, y_2, 0)) \geq \max\{2\log|x_2 - x_1|, 2\log|y_2 - y_1|\}$.

*Proof.* See Lemma 3.5.2 in [1]. $\qquad\square$

We now obtain a lower bound for $d_S$ in terms of the Euclidean distance:

$$d_{\mathbb{R}^2}((x_1, y_1), (x_2, y_2)) \leq \sqrt{2}\max\{|x_2 - x_1|, |y_2 - y_1|\} \qquad (3.4)$$

$$\log d_{\mathbb{R}^2}((x_1, y_1), (x_2, y_2)) \leq \max\{\log|x_2 - x_1|, \log|y_2 - y_1|\} + \log\frac{1}{\sqrt{2}} \quad (3.5)$$

hence by Lemma 3.7.7

$$d_S((x_1, y_1, 0), (x_2, y_2, 0)) \geq 2\log d_{\mathbb{R}^2}((x_1, y_1, 0), (x_2, y_2, 0)) + \log 2 \quad (3.6)$$

This estimate bounds the distortion of $\mathbb{R}^2$ in $Sol$ and our next goal will be to quasify it to apply to the discrete case of $P_2$. Before we do that however, we provide an explicit embedding of the non-nilpotent group $P_{2,\phi}$ as a uniform lattice in $Sol$ following [1].

**The Embedding of $P_{2,\phi}$ in $Sol$**

First, we embed $\mathbb{Z}^2$ as the standard integer lattice in $\mathbb{R}^2$. Let the eigenvalues of $\phi$ be $0 < \lambda^{-1} < 1 < \lambda$, and let $\mathbf{v}_-$ and $\mathbf{v}_+$ be eigenvectors for the corresponding eigenvalues. Let $f : \mathbb{R}^2 \to \mathbb{R}^2$ be the linear map defined by $\mathbf{v}_- \mapsto e_1$ and $\mathbf{v}_+ \mapsto e_2$. Now, we define a map $\psi : P_{2,\phi} \to Sol$ by $\psi(e_1^k e_2^l t^q) = (f(k, l), q\log(\lambda))$, where $e_1, e_2$ and $t$ are the standard generators of $P_{2,\phi}$ described in Section 3.3, and $k, l, q \in \mathbb{Z}$. It is not hard to show that the map $\psi$ is a monomorphism with a discrete image and that $\psi(P_{2,\phi})\backslash Sol$ is compact.

**Bounding the Growth of $(P_{2,\phi}, \langle t \rangle)$**

By uniformity of the embedding, we have the following quasi-isometry inequality:

$$d_{P_{2,\phi}}((x_1, y_1, 0), (x_2, y_2, 0)) \geq Cd_S((x_1, y_1, 0), (x_2, y_2, 0)). \quad (3.7)$$

Now, passing from the continuous metric on $\mathbb{R}^2$ to the word metric on $\mathbb{Z}^2$ we have

$$d_{\mathbb{R}^2}((x_1, y_1, 0), (x_2, y_2, 0)) \geq \frac{1}{\sqrt{2}}d_{\mathbb{Z}^2}((x_1, y_1, 0), (x_2, y_2, 0)) \quad (3.8)$$

and finally putting all these estimates together we get:

$$d_{P_{2,\phi}}((x_1, y_1, 0), (x_2, y_2, 0)) \geq C \log d_{\mathbb{Z}^2}((x_1, y_1, 0), (x_2, y_2, 0)) + A \qquad (3.9)$$

**Proposition 3.7.8.** *Let $f_H(r)$ be the coset growth of $H = \langle t \rangle$ in $P_{2,\phi}$, then $f_H(r) \sim \alpha^r$.*

*Proof.* Given $r > 0$, Lemma 3.5.2 shows that there exists $N \sim 2\left(\pi r^2\right)^2$ such that any $N$ distinct cosets of $H$, say $a_1 H, ..., a_N H$, contains a pair $aH, bH$ with

$d_{\mathbb{Z}^2}(\phi^{-k}(a), \phi^{-k}(b)) > r$. Now, we have $d_{P_{2,\phi}}(at^k, bt^k) = d_{P_{2,\phi}}(\phi^{-k}(a), \phi^{-k}(b))$, which combined with (3.9) gives $d_{P_{2,\phi}}(at^k, bt^k) > C \log r + A$. Arguing as in the proof of Lemma 3.6.1, we conclude that $d_{P_{2,\phi}}(aH, bH) \geq \frac{1}{2}(C \log r + A)$. Now, setting $R = \frac{1}{4}(C \log r + A) - 1$, we conclude that the open ball of radius $R$ in $P_{2,\phi}/H$ has at most $N \sim B\alpha^r$ elements, where $B$ and $\alpha$ depend only on $A$ and $C$. $\qquad \square$

# Chapter 4

# On the Quasiconvex Subgroups of $F_m \times \mathbb{Z}^n$

## 4.1 Introduction and Basic Notions

The motivation for the work in this chapter comes from the following remark found in the introductory section of [13]: "...it is currently unknown whether a quasiconvex subgroup of a CAT(0) group is itself CAT(0)." At the time of first being acquainted with the problem, it seemed to me almost absurd that the answer to such a basic question in the theory of CAT(0) groups was not yet known. By comparison, the corresponding statement in the theory of hyperbolic groups has long been known to be true, see Lemma 3.1.3.

Let us begin by recalling the definition of quasiconvexity from Section 2.2.1:

**Definition 4.1.1.** Let $X$ be a (uniquely) geodesic metric space and let $Y \subseteq X$ be a subspace. We shall say that $Y$ is $\nu$-quasiconvex if there exists $\nu > 0$ such that $[x, y] \subseteq \mathcal{N}_\nu(Y)$, for all $x, y \in Y$.

Quasiconvexity is a generalization, or more precisely, a "quasification" of the notion of convexity in a geodesic space. Recall that a subspace is called convex if geodesics joining points in that subspace are completely contained in it. By contrast, the notion of quasiconvexity allows for some wiggle room: a geodesic joining points in $Y$ need not be contained in $Y$ itself but is rather allowed to travel in a fixed bounded neighborhood of $Y$ instead.

**Definition 4.1.2.** Let $G$ be a CAT(0) group acting geometrically on the CAT(0) space $X$, let $H \subseteq G$ be a subgroup, and let $x_0 \in X$ be a basepoint. The subgroup $H$ is called $\nu$-quasiconvex if the group orbit $Hx_0$ is a $\nu$-quasiconvex subspace of $X$.

Often, we make no mention of the quasiconvexity constant $\nu$ and simply say that $H$ is a *quasiconvex* subgroup. While the value of the constant $\nu$ in Definition 4.1.2 may depend on the choice of the basepoint $x_0$, whether $H$ is quasiconvex or not does not depend on this choice.

Let us compare Definition 4.1.2 with the notion of quasiconvexity in hyperbolic groups given in Definition 3.1.1. The most obvious difference between the two is that Definition 3.1.1 is an intrinsic definition, it makes no mention of anything outside the group $G$, the subgroup $H$, and a condition on $H$ in terms of the word metric on $G$, while Definition 4.1.2 defines quasiconvexity in terms of a group orbit in a given CAT(0) space. This difference is the reason why the notion of quasiconvexity in hyperbolic groups is more robust than its CAT(0) counterpart: as the following example found in [13] shows, considering two actions of $G = F_2 \times \mathbb{Z}$ on the same CAT(0) space, we can arrange a subgroup of $G$ to be quasiconvex with respect to one action but not the other.

**Example 4.1.1.** We consider the group $G = F_2 \times \mathbb{Z} = \langle a, b \rangle \times \langle t \rangle$. Let $X$ be the universal cover of the presentation 2-complex of $G$ with the usual action of $G$. If all the edges in $X$ have length equal to 1, it is easy to see that the subgroup $H = \langle a, b \rangle$ is $\nu$-quasiconvex for any $\nu > \frac{1}{2}$. Note that the subgroup $K = \langle a, bt \rangle$ is not quasiconvex, otherwise the subgroup $H \cap K$ would also be quasiconvex, as the intersection of two quasiconvex subgroups is quasiconvex. However, the subgroup $H \cap K$ is the subgroup of $H$ which consists of elements for which the sum of the powers of the generator $b$ in a reduced word equals 0. Since quasiconvex $\Rightarrow$ finitely generated by Lemma 3.1.3, $H \cap K$ is not quasiconvex as it is not finitely generated.

Now, consider the automorphism $\phi$ of $G$ defined by $a \mapsto a$, $b \mapsto bt$, $t \mapsto t$, which sends $H$ to $K$. If we compose the action $\rho : G \to Isom(X)$ of $G$ on $X$ described above with $\phi$, we obtain a new action $\rho \circ \phi : G \to Isom(X)$, such that $\rho$ restricted to $K$ is the same as $\rho \circ \phi$ restricted to $H$. Thus, $H$ is not quasiconvex with respect to this new action.

In order to show that a quasiconvex subgroup $H$ of the CAT(0) group $G$ is CAT(0), we need to exhibit a CAT(0) space $Y$ and a geometric action of $H$ on $Y$. The reader may at this point rightfully ask the question, why can we not take the convex hull of the $H$-orbit of some point $x_0 \in X$? The convex hull is a convex, $H$-invariant subspace of a CAT(0) space, and the action of $H$ on it is proper, as the action of $G$ on $X$ is proper. The problem is this: while the action of $G$ on $X$ was assumed to be cocompact, it is not at all obvious, and perhaps not even true in general, that quasiconvexity of $H$ in $G$ should imply cocompactness of the action of $H$ on $conv(Hx_0)$. Of course, even if one were to find a counterexample, that is an example where the induced action of $H$ on $conv(Hx_0)$ is not cocompact, this would not necessarily mean that $H$ is not a CAT(0) group because there is still the possibility that $H$ may act geometrically on some other CAT(0) space. As far as we know, the general question of whether quasiconvexity implies CAT(0) is still wide open. In this chapter, we show that what *should* be true is indeed true in one special case, namely, we prove the following:

**Theorem**: Let $H$ be a quasiconvex subgroup of $G = F_m \times \mathbb{Z}^n$, and let $X$ be the product of the regular $2m$-valent tree with $\mathbb{R}^n$ with the usual action of $G$. Then the action of $H$ on the convex hull of any orbit $Hx_0$ is cocompact.

If $G$ is a CAT(0) group acting geometrically on the CAT(0) space $X$, and there exists a closed convex $H$-invariant subset of $X$ on which $H$ acts cocompactly, $H$ is called *convex*.

With this terminology our theorem becomes: *Any quasiconvex subgroup of $G = F_m \times \mathbb{Z}^n$ is convex with respect to the usual action of $G$ on $Tree \times \mathbb{R}^n$.*

In the course of proving this result, we introduce a technique for analyzing the convex hull in certain CAT(0) spaces and also prove that quasiconvex subgroups of $F_m \times \mathbb{Z}^n$ are virtually of the form $A \times B$ where $A \le F_m$ and $B \le \mathbb{Z}^n$, which may be regarded as a structure result of sorts.

## 4.2 Convex Hulls and Quasiconvex Subgroups

As we mentioned in the previous section, every convex subspace of a CAT(0) space is obviously itself a CAT(0) space with the induced metric. An idea which dates back to Minkowski and Brunn, and which we independently rediscovered, is to construct $conv(Y)$ by means of a sequential process as follows: For $S \subseteq X$, we define $conv^1(S)$ to be the union of all geodesic segments having both endpoints in $S$, or symbolically $conv^1(S) = \bigcup_{s_1, s_2 \in S} [s_1, s_2]$. Now, we set $conv^0(Y) = Y$ and define recursively $conv^i(Y) = conv^1(conv^{i-1}(Y))$. This process of "convexification" results in an ascending sequence of subsets of $X$: $Y = conv^0(Y) \subseteq conv^1(Y) \subseteq \ldots \subseteq conv^i(Y) \subseteq \ldots \subseteq conv(Y) \subseteq X$, each of which gets closer to the convex hull of $Y$ in the following sense:

**Lemma 4.2.1.** *Let $X$ be a uniquely geodesic space and $Y$ a subspace, then $conv(Y) = \bigcup_{i=0}^{i=\infty} conv^i(Y)$.*

*Proof.* Obviously, $conv^i(Y) \subseteq conv(Y)$ for all $i$. It is also clear that $\bigcup_{i=0}^{i=\infty} conv^i(Y)$ is convex, hence it equals $conv(Y)$. $\square$

Lemma 4.2.1 was the starting point for our investigation of convex hulls. Its proof is not difficult and we realized that the result was already attributed to Hermann Brunn in the case when $X$ is a vector space only after we had proved the results on convexity in polygonal complexes below.

**Remark 4.2.2.** A set $Y$ is $\nu$-quasiconvex if and only if $conv^1(Y) \subseteq \mathcal{N}_\nu(Y)$.

We now relate the foregoing discussion on convexity with non-positive curvature. In CAT(0) spaces, as a consequence of the convexity of the metric, one has control over the growth of the sizes of the sets $conv^i(Y)$ as the following result shows:

**Lemma 4.2.3.** *Let $X$ be a CAT(0) space and let $Y \subseteq X$ be a $\nu$-quasiconvex subset of $X$. Then, $conv^i(Y) \subseteq N_{i\nu}(Y)$.*

*Proof.* In this proof we assume that all geodesics are parametrized proportional to arc length and we proceed by induction on $i$. The starting step $i = 1$ is handled by Remark 4.2.2. For the inductive step, suppose that $conv^{i-1}(Y) \subseteq N_{(i-1)\nu}(Y)$ and let $x \in conv^i(Y)$. Then, $x \in [x_1, x_2]$, where $x_1, x_2 \in conv^{i-1}(Y) \subseteq N_{(i-1)\nu}(Y)$. Let $x_1', x_2' \in Y$ be such that $d(x_j, x_j') < (i-1)\nu$, for $j = 1, 2$. By convexity of the CAT(0) metric, see Proposition 2.2.4, $d([x_1, x_2](t), [x_1', x_2'](t)) \leq (1-t)d(x_1, x_1') + td(x_2, x_2') < (i-1)\nu$. This shows that $d(x, conv^1(Y)) < (i-1)\nu$ and thus we conclude that $x \in N_{i\nu}(Y)$, as desired. $\square$

The number $k = inf\{i : conv^i(Y) = conv(Y)\}$ is called the Brunn number. Brunn gave a lower and an upper bound for $k$ in finite-dimensional vector spaces, see [4]. In view of Lemma 4.2.3, a uniform bound for the Brunn number over all subsets of a CAT(0) space has important implications from the point of view of cocompactness of group actions. Before we proceed to give a bound for the Brunn number for Euclidean spaces, we make the following important observation: In analyzing the convex hull of an arbitrary set using Lemma 4.2.1, it suffices to only consider finite sets, which are substantially easier to work with.

**Lemma 4.2.4.** *Let $X$ be a geodesic space, let $Y \subseteq X$, and let $i \in \mathbb{N}^+$. If $conv^i(S) = conv(S)$ for every finite subset $S \subseteq Y$, then $conv^i(Y)$ is convex.*

*Proof.* Let $x, y \in conv^i(Y)$. Then, we can find $a_1, a_2, b_1, b_2 \in conv^{i-1}(Y)$ such that $x \in [a_1, a_2]$ and $y \in [b_1, b_2]$. Similarly, we can find $c_1, c_2, d_1, d_2 \in conv^{i-2}(Y)$ such that $a_1 \in [c_1, c_2]$ and $a_2 \in [d_1, d_2]$, etc. Proceeding recursively, we see that we can find points $x_1, \ldots, x_m \in Y$ such that $x, y \in conv^i(\{x_1, \ldots, x_m\}) = conv(\{x_1, \ldots, x_m\})$. Hence, $[x, y] \subseteq conv^i(\{x_1, \ldots, x_m\}) \subseteq conv^i(Y)$. $\qquad\square$

As we were unable to procure Brunn's original paper, we present our own proof of the intuitively obvious fact that the Brunn number of any subset of $\mathbb{R}^n$ is less or equal to $n$. Our proof uses Caratheodory's theorem which we recall below.

**Convex Hulls in Euclidean Spaces**

**Theorem 4.2.5.** *(Caratheodory) If $E$ is a vector space of dimension $d$, then, for every subset $X$ of $E$, every element in the convex hull $conv(X)$ is an affine convex combination of $d + 1$ elements of $X$.*

*Proof.* See Proposition 5.2.3 in [21]. $\qquad\square$

**Lemma 4.2.6.** *For any finite set $S \subseteq \mathbb{R}^n$, $conv^n(S) = conv(S)$.*

*Proof.* By Caratheodory's theorem, $conv(S) = \bigcup conv(\{s_1, ..., s_{n+1}\})$, where the union is taken over all $s_1, ..., s_{n+1} \in S$. Therefore, it suffices to show that for any set $\{s_1, ..., s_{n+1}\} \subseteq \mathbb{R}^n$, $conv^n(\{s_1, ..., s_{n+1}\}) = conv(\{s_1, ..., s_{n+1}\})$. Consider the points $e_1, ..., e_{n+1} \in \mathbb{R}^{n+1}$. Their convex hull is the standard $n$-simplex $\Delta_n$ in $\mathbb{R}^{n+1}$. Suppose $conv^{i-1}(\{e_1, ..., e_{n+1}\})$ contains all the $(i-1)$-faces of $\Delta_n$, then $conv^i(\{e_1, ..., e_{n+1}\})$ contains all joins of the form $join\{F, e_j\}$, $1 \le j \le n+1$, where $F$ is an $(i-1)$-face. But all of the $i$-faces are joins of this form. By induction, $conv^i\{e_1, ..., e_{n+1}\}$ contains all the $i$-faces. Hence, $conv^n(\{e_1, ..., e_{n+1}\})$ contains and therefore equals $\Delta_n$. Now, let $\phi$ be the affine map which sends $e_i$ to $s_i$. This map sends lines to lines, therefore $\phi(conv^i(\{e_1, ..., e_{n+1}\})) \subseteq conv^i(\{s_1, ..., s_{n+1}\})$. Now, $\phi(conv^n(\{e_1, ..., e_{n+1}\}))$ is convex, contains $\{s_1, ..., s_{n+1}\}$, and is contained in $conv^n(\{s_1, ..., s_{n+1}\})$. Therefore, $conv(\{s_1, ..., s_{n+1}\}) = conv^n(\{s_1, ..., s_{n+1}\})$, as desired. $\qquad\square$

Combining Lemma 4.2.4 and Lemma 4.2.6, we obtain the desired bound on the Brunn number:

**Corollary 4.2.7.** *For any subset $Y \subseteq \mathbb{R}^n$, $conv^n(Y) = conv(Y)$. Therefore, $k \le n$ for any subset of $\mathbb{R}^n$.*

As a straightforward application of Corollary 4.2.7 we have:

**Corollary 4.2.8.** *In any uniquely geodesic Hilbert geometry, the Brunn number of any subset is bounded above by the dimension of the underlying Euclidean space.*

*Proof.* As we saw in Example 2.2.3, each affine segment in the underlying Euclidean space is a geodesic for the Hilbert metric. Since the Hilbert geometry we consider is uniquely geodesic, Corollary 4.2.7 immediately applies. $\qquad\square$

Unfortunately, obtaining a bound on the Brunn number in an arbitrary CAT(0) space in the absence of any restrictions on the subsets in consideration is impossible. We are, however, able to bound the Brunn number in certain "planar" piecewise Euclidean CAT(0) polygonal complexes.

**Convex Hulls in CAT(0) Planes**

**Definition 4.2.1.** A *CAT(0) plane* is a simply connected piecewise Euclidean polygonal complex $X$ with $Shapes(X)$ finite, such that $Lk(v)$ is isometric to a circle of length $\ge 2\pi$ for every vertex $v \in X$.

We immediately note that in view of Theorem 2.2.13, a CAT(0) plane is a CAT(0) metric space. We prove that the Brunn number of any subset of a CAT(0) plane is $\leq 2$. The proof of this intuitively "obvious" result turned out to be surprisingly technical. Our proof employs a local-to-global technique which makes use of the Cartan-Hadamard theorem. The idea is that under mild hypotheses, convexity on the small scale implies global convexity. Before we are able to state and prove the promised results, we need to transpose some familiar definitions to the "small scale":

**Definition 4.2.2.** Let $(X, d)$ be a metric space.

1. The metric $d$ is said to be *locally convex* if every point in $X$ has a neighborhood in which the induced metric is convex.

2. The metric space is said to be *locally CAT(0)* if every point in $X$ has a convex neighborhood $\mathcal{U}$ with the property that $(\mathcal{U}, d)$ is a CAT(0) metric space.

3. Let $f : X \to Y$ be a map between two metric spaces. We shall say that $f$ is a *local isometry* if every point $x \in X$ has a neighborhood $\mathcal{U}$, such that $f$ restricted to $\mathcal{U}$ is an isometric embedding.

**Theorem 4.2.9.** *(Cartan-Hadamard) Let $(X, d)$ be a complete connected metric space. If the metric on $X$ is locally convex, the the induced length metric on the universal covering $\widetilde{X}$ is convex. In particular, there is a unique geodesic segment joining each pair of points in $\widetilde{X}$. Further, if $X$ is a locally CAT(0) space, then $\widetilde{X}$ with the induced length metric is a CAT(0) space and the covering map $p : \widetilde{X} \to X$ is a local isometry.*

*Proof.* See Theorem 4.1 in [3]. $\qquad\qquad\qquad\qquad\qquad\qquad\qquad\qquad\qquad\qquad\square$

In view of the Cartan-Hadamard theorem above, we shall say that a subset $Y$ of the CAT(0) space $X$ is *locally convex* if every point of $Y$ has a convex neighborhood $\mathcal{U} \subseteq X$ such that $\mathcal{U} \cap Y$ is convex.

**Proposition 4.2.10.** *Let $X$ be a CAT(0) space. Then, $conv^n(Y) = conv(Y)$ for all $Y \subseteq X$ if and only if $conv^n(S)$ is locally convex for every finite subset $S \subseteq X$.*

*Proof.* First, we show that for any compact subset $S \subseteq X$, $conv^i(S)$ is compact for every $i$. The proof is by induction on $i$, the case $i = 0$ being trivial. Note that we have an obvious map $\varphi : conv^{i-1}(S) \times conv^{i-1}(S) \times I \to conv^i(S)$ given by $\varphi(x, y, t) = [x, y](t)$. Since in any CAT(0) space geodesics vary continuously with endpoints, the map $\varphi$ is a continuous surjection which maps the compact set $conv^{i-1}(S) \times conv^{i-1}(S) \times I$ to $conv^i(S)$ thus proving our claim.

Now, suppose that $conv^n(S)$ is locally convex for every finite subset $S \subseteq X$. Since $conv^n(S)$ is a compact and therefore complete, connected, and locally convex subset of a CAT(0) space, Theorem 4.2.9 tells us that the universal cover $\widetilde{conv^n(S)}$ endowed with the length metric is a CAT(0) space, and that the covering map $p : \widetilde{conv^n(S)} \to conv^n(S)$ is a local isometry. Let $x, y \in conv^n(S)$ and choose any $\tilde{x} \in p^{-1}(x)$, and $\tilde{y} \in p^{-1}(y)$. Let $\alpha(t)$ be the unique geodesic in $\widetilde{conv^n(S)}$ joining $\tilde{x}$ to $\tilde{y}$. Then, because $conv^n(S)$ is compact, a simple argument using the Lebesgue covering lemma shows that $p \circ \alpha(t)$ is a local geodesic in $X$ joining $x$ to $y$. Since in a CAT(0) space any local geodesic is a geodesic, see Proposition 1.4 in [3], we see that $p \circ \alpha$ is the geodesic in $X$ joining $x$ to $y$. But the image of $p \circ \alpha$ is contained in $conv^n(S)$. This shows that $conv^n(S)$ is convex, and Lemma 4.2.4 yields the desired conclusion. $\square$

In the course of proving Proposition 4.2.13 we will encounter the phenomenon of *bifurcating* geodesics. We shall call geodesics which coincide up to a point of *bifurcation* or divergence, bifurcating geodesics. The key fact we shall need is that in a CAT(0) plane, at the point of bifurcation, the geodesic which splits can be extended in an infinite number of ways. To prove this result, we need the following reformulation of the CAT(0) condition:

**Lemma 4.2.11.** *Let $X$ be a geodesic metric space. Then, $X$ is a CAT(0) space if and only if the Alexandrov angle (see Definition 2.2.4) between the sides of any geodesic triangle in $X$ with distinct vertices is no greater than the angle between the corresponding sides of its comparison triangle in $\mathbb{E}^2$.*

*Proof.* See Proposition 1.7(4) in [3]. □

Throughout the proof of Lemma 4.2.12, we employ notation from Section 2.2.2.

**Lemma 4.2.12.** *Let $X$ be a CAT(0) plane and let $a_1, a_2, b, q \in X$ be such that $[a_1, b] \cap [a_2, b] = [q, b]$. Then, for any $c \in [a_1, a_2]$, the concatenation of the geodesic segment $[c, q]$ and $[q, b]$ is a geodesic segment.*

*Proof.* The proof is essentially the observation that for any $c \in [a_1, a_2]$, the distance in $Lk(q)$ between the directions determined by $[c, q]$ and $[q, b]$ is at least $\pi$, hence the Alexandrov angle $\angle_q(c, b) = \pi$. By Lemma 4.2.11, $\angle_q(c, b) \leq \overline{\angle}_{\bar{q}}(\bar{c}, \bar{b})$, hence $\overline{\angle}_{\bar{q}}(\bar{c}, \bar{b}) = \pi$. Now, we conclude that $d(\bar{c}, \bar{b}) = d(\bar{c}, \bar{q}) + d(\bar{q}, \bar{b})$. Therefore, $d(c, b) = d(c, q) + d(q, b)$ thus showing that the concatenation of $[c, q]$ and $[q, b]$ is indeed a geodesic segment. □

We are now ready to show that the Brunn number for any subset of a CAT(0) plane is at most equal to 2.

**Proposition 4.2.13.** *Let $X$ be a CAT(0) plane. Then the Brunn number of any subset $Y \subseteq X$ is at most 2.*

*Proof.* In view of Proposition 4.2.10, we only need to show that for every finite subset $S \subseteq X$, $conv^2(S)$ is locally convex at every point $p \in conv^2(S)$. Suppose to the contrary that $conv^2(S)$ is not locally convex at some $p \in conv^2(S)$. Then, given any convex neighborhood $\mathcal{U}$ of $p$ in $X$, since $conv^2(S)$ is compact, we can find a geodesic $\gamma : [0, 1] \to \mathcal{U}$ such that $\gamma(0), \gamma(1) \in conv^2(S)$ and $\gamma(t) \notin conv^2(S)$ for all $0 < t < 1$. Let $\gamma(0) = x_0$ and note that we may have $x_0 = p$. Then, there exist points $x_1, x_2, y_1, y_2 \in S$ such that $x_0 \in [[x_1, y_1](t_1), [x_2, y_2](t_2)]$ for some $t_1, t_2 \geq 0$. Consider the 1-parameter family of geodesics $t' \mapsto [[x_1, y_1](t'), [x_2, y_2](t_2)]$. If $x_0 \in [[x_1, y_1](t'), [x_2, y_2](t_2)]$ for all $t' \geq t_1$ (or all $t' \leq t_1$), then $x_0$ lies on a geodesic having one endpoint $s$ in $S$. In this case we consider the family $t' \mapsto [s, [x_2, y_2](t')]$.

If $x_0 \in [s, [x_2, y_2](t')]$ for all $t' \geq t_2$ (or all $t' \leq t_2$), then $x_0 \in conv^1(S)$, which is a contradiction. Therefore, without loss of generality, we may assume that $x_0 \notin [[x_1, y_1](t'), [x_2, y_2](t_2)]$ for some values of $t'$ both greater than and less than $t_1$. As $X$ has no free edges, $\gamma$ may be extended to $\gamma : [-\epsilon, 1] \to X$ such that this extension of $\gamma$ crosses $[[x_1, y_1](t_1), [x_2, y_2](t_2)]$ at $x_0$. Because $\gamma(t) \notin conv^2(S)$ for $t > 0$ and because geodesics vary continuously with endpoints, we can find $t'_0 < t_1 < t''_0$ and $-\epsilon < -\epsilon' < 0$ such that $[[x_1, y_1](t'_0), [x_2, y_2](t_2)]$ and $[[x_1, y_1](t''_0), [x_2, y_2](t_2)]$ intersect at $\gamma(-\epsilon')$, and such that neither of these geodesic segments passes through $x_0$. Because of uniqueness of geodesics in $X$, the segments $[[x_1, y_1](t'_0), [x_2, y_2](t_2)]$ and $[[x_1, y_1](t''_0), [x_2, y_2](t_2)]$ must coincide up to some point $q$ which is their point of bifurcation. Then, Lemma 4.2.12 shows that for every $t_1 < t < t_2$, the concatenation of the geodesic segments $[[x_2, y_2](t_2), q]$ and $[q, [x_1, y_1](t)]$ is a geodesic. Again by uniqueness of geodesics, this concatenation must be the geodesic segment $[[x_1, y_1](t_1), [x_2, y_2](t_2)]$, and we conclude that $x_0 \in [[x_1, y_1](t'_0), [x_2, y_2](t_2)]$ which is a contradiction. $\square$

## 4.3  Free×Free-abelian Groups

Throughout this section, $G$ will be the group $F_m \times \mathbb{Z}^n$ and $X$ will stand for the product of the regular $2m$-valent tree $T_{2m}$ and $\mathbb{R}^n$. The action of $G$ on $X$ is the product action where $F_{2m}$ acts as the group of deck transformations on the universal cover of the m-rose, $T_{2m}$, and $\mathbb{Z}^n$ acts by translation on $\mathbb{R}^n$.

**Lemma 4.3.1.** *Let $H = \langle f_1 z_1, ..., f_s z_s \rangle$, $f_i \in F_m$, $z_i \in \mathbb{Z}^n$ be a quasiconvex subgroup of $G$ such that not all of the $f_i$ have the same axis of translation in $T_{2m}$. Then, there exist positive integers $k_1, ..., k_s$ such that $H$ contains the subgroup $A = \langle z_1^{k_1}, ..., z_s^{k_s} \rangle$.*

*Proof.* Let $1 \leq i \leq s$, let $j$ be such that $f_i$ and $f_j$ have different axes of translation, and let $l$ be a positive integer. Find an axis of translation for $f_i z_i$ whose projection to the Euclidean factor passes through $0 \in \mathbb{R}^n$. This can always be done by translating the Euclidean component of any given axis for $f_i z_i$. Let $x_0$ be a point on the chosen axis of translation for $f_i z_i$ such that $pr_{\mathbb{R}^n}(x_0) = 0$. Consider the sequences of points $(f_i z_i)^l x_0$ and $(f_j z_j)(f_i z_i)^l x_0$. Because $f_i$ and $f_j$ have different axes, there is a vertex $v$ in $T_{2m}$ such that for every $l$, the geodesic segment $\big[(f_i z_i)^l x_0, (f_j z_j)(f_i z_i)^l x_0\big]$ passes through the flat $\{v\} \times \mathbb{R}^n$. Let $y_l$ denote the point of intersection of $\{v\} \times \mathbb{R}^n$ and $\big[(f_i z_i)^l x_0, (f_j z_j)(f_i z_i)^l x_0\big]$. The orthogonal projection of this geodesic segment to the flat $\{v\} \times \mathbb{R}^n \cong \mathbb{R}^n$ is the geodesic segment between $z_i^l$ and $z_i^l + z_j$. Since the geodesic $\big[(f_i z_i)^l x_0, (f_j z_j)(f_i z_i)^l x_0\big]$ intersects its orthogonal projection in the point $y_l$, we have $d(y_l, (v, z_i^l)) \leq \|z_j\|$. The orbit $Hx_0$ is quasiconvex, hence there is $\nu > 0$ and $h_l \in H$ such that $d(h_l x_0, y_l) < \nu$. Then, $d(h_l x_0, z_i^l x_0) \leq d(h_l x_0, y_l) + d(y_l, z_i^l(v,0)) + d(z_i^l(v,0), z_i^l x_0) < \nu + \|z_j\| + d(x_0, (v,0))$. If $\tau = \nu + \|z_j\| + d(x_0, (v,0))$, then $B_\tau(h_l x_0) \cap B_\tau(z_i^l x_0) \neq \emptyset$, or $B_\tau(h_l^{-1} z_i^l x_0) \cap B_\tau(x_0) \neq \emptyset$, for all $l$. Because the action of $G$ is proper, $h_l^{-1} z_i^l = g \in G$ for infinitely many values of $l$. Then, for some $k, l$ we have $z_i^{l-k} = h_l h_k^{-1} \in H$. Setting $k_i = l - k$, we obtain $z_i^{k_i} \in H$. $\qquad \square$

Let $p : T_{2m} \times \mathbb{Z}^n \to T_{2m}$ be the projection onto the first factor, and let $V$ denote the real span of the vectors $z_1^{k_1}, ..., z_k^{k_s}$. We then, have the following:

**Lemma 4.3.2.** *With the same notation as in Lemma 4.3.1, the convex hull of $Hx_0$ equals $conv(p(Hx_0)) \times V$.*

*Proof.* First, we note that the projection maps $p, pr_{\mathbb{R}^n}$ commute with the operation of forming the convex hull. That is, $p(conv(Hx_0)) = conv(p(Hx_0))$, and similarly for $pr_{\mathbb{R}^n}$. Let us show this for the projection map $p$. We begin by making the observation that $p(conv^1(S)) = conv^1(p(S))$ for any set $S$, since $p$ maps the geodesic segment connecting two points to the geodesic segment connecting their images. Therefore, we have $p(conv(Hx_0)) = p\left(\bigcup_i conv^i(Hx_0)\right) = \bigcup_i p(conv^i(Hx_0)) = \bigcup_i conv^i(p(Hx_0)) = conv(p(Hx_0))$.

Now, we proceed with the proof of the lemma.

'$\subseteq$': Without loss of generality, we may assume that $pr_{\mathbb{R}^n}(x_0) = 0$. Clearly, $conv(Hx_0) \subseteq p(conv(Hx_0)) \times pr_{\mathbb{R}^n}(conv(Hx_0))$, which after commuting the projection maps past $conv$ gives us the desired inclusion.

'$\supseteq$': Let $x \in conv(p(Hx_0)) \times V$. Let $y \in conv(Hx_0)$ be such that $p(y) = p(x)$. Note that because $H$ contains powers of the Euclidean translations $z_1, ..., z_k$, the projection of the convex hull of the orbit $Hx_0$ to the Euclidean factor will equal $V$. Also, $conv(Hx_0) \supseteq V \cdot y$, as $conv(Hx_0)$ is stable under the action of $V$ by translations on the second factor. Hence, we can write $x = w \cdot y$, for some $w \in V$, so that $x \in conv(Hx_0)$. $\qquad\square$

**Lemma 4.3.3.** *Let $H$ be as in Lemma 4.3.1. Then, the group $H$ acts cocompactly on its convex hull. In particular, the Brunn number of the orbit $Hx_0$ is bounded above by $1 + \dim(V)$.*

*Proof.* Note that we can write $p(conv(Hx_0)) = conv(p(Hx_0))$ as a union of biinfinite geodesic rays $\gamma$, such that any point on $\gamma$ lies between two points in $p(Hx_0)$. Then, $conv(Hx_0) = \bigcup_\gamma \gamma \times V$, and $\gamma \times V \cong \mathbb{R}^{1+\dim(V)}$. Note that $\gamma \times V$ contains the lattices $p(hx_0) \times \mathbb{Z}\text{-}span \left\langle z_1^{k_1}, ..., z_s^{k_s} \right\rangle$, where $h \in H$. Because of the assumption that any point on $\gamma$ lies between two points $p(h_1), p(h_2) \in p(Hx_0)$, the convex hull of these lattices is all of $\gamma \times V$, and by Corollary 4.2.7, $conv^{1+\dim(V)} \left( \bigcup(p(hx_0) \times \mathbb{Z}\text{-}span \left\langle z_1^{k_1}, ..., z_s^{k_s} \right\rangle \right) = \gamma \times V$. Finally, $conv^{1+\dim(V)}(Hx_0) \supseteq conv^{1+\dim(V)} \left( \bigcup(p(hx_0) \times \mathbb{Z}\text{-}span \left\langle z_1^{k_1}, ..., z_s^{k_s} \right\rangle) \right) = conv(Hx_0)$. $\qquad\square$

Combining Lemmas 4.3.1-4.3.3, we obtain:

**Theorem 4.3.4.** *Any quasiconvex subgroup of $F_m \times \mathbb{Z}^n$ acts cocompactly on the convex hull of any of its orbits.*

*Proof.* Lemmas 4.3.1-4.3.3 take care of the case when for each $i$ there is $j$ such that such that $f_i$ and $f_j$ have different axes of translation. If $H = \left\langle f^{k_1} z_1, ..., f^{k_s} z_s \right\rangle$, $f \in F_m$, then the orbit $Hx_0$ is contained in a single flat $a_f \times V$ isometric to $\mathbb{R}^{1+\dim(V)}$, where $a_f$ is a common axis for all $f^{k_i}$, and $x_0$ is on a common axis for all the $f_i z_i$. Hence, $conv(Hx_0) = conv^{1+\dim(V)}(Hx_0)$, which shows cocompactness of the action of $H$. In either of the cases $H = \langle f_1, ..., f_s \rangle$ or $H = \langle z_1, ..., z_s \rangle$, the conclusion is again trivially true. In the former case $conv(Hx_0) = conv^1(Hx_0)$, while in the latter $conv(Hx_0) = conv^s(Hx_0)$. $\qquad\square$

In the course of proving the theorem, we have the essential ingredients for the following corollary:

**Corollary 4.3.5.** *If $H$ is a quasiconvex subgroup of $F_m \times \mathbb{Z}^n$, then $H$ is virtually of the form $A \times B$, where $A \leq F_m$ and $B \leq \mathbb{Z}^n$.*

*Proof.* Let $H$ be as in Lemma 4.3.1. The proof of Lemma 4.3.1 shows that for $g = fz \in H$, there exists $s$ such that $z^s \in H$, and hence $f^s \in H$. Let $A = \mathbb{Z}^n \cap H$, $F = F_m \cap H$. Then, $g^s \in AF$. On the other hand, $[H, H] \subseteq F$, and also $AF$ is normal in $H$. Hence, we see that $H/AF$ is a finitely generated, torsion, abelian group, and is therefore finite. If $H = \left\langle f^{k_1} z_1, ..., f^{k_s} z_s \right\rangle$, or $H = \langle z_1, ..., z_s \rangle$, then $H$ is already free abelian. $\qquad\square$